# On Universality and Training in Binary Hypothesis Testing

Michael Bell and Yuval Kochman

## Abstract

The classical binary hypothesis testing problem is revisited. We notice that when one of the hypotheses is composite, there is an inherent difficulty in defining an optimality criterion that is both informative and well-justified. For testing in the simple normal location problem (that is, testing for the mean of multivariate Gaussians), we overcome the difficulty as follows. In this problem there exists a natural "hardness" order between parameters as for different parameters the error-probabilities curves (when the parameter is known) are either identical, or one dominates the other. We can thus define minimax performance as the worst-case among parameters which are below some hardness level. Fortunately, there exists a universal minimax test, in the sense that it is minimax for all hardness levels simultaneously. Under this criterion we also find the optimal test for composite hypothesis testing with training data. This criterion extends to the wide class of local asymptotic normal models, in an asymptotic sense where the approximation of the error probabilities is additive. Since we have the asymptotically optimal tests for composite hypothesis testing with and without training data, we quantify the loss of universality and gain of training data for these models.

The authors are with the Rachel and Salim Benin school of computer science and engineering, Hebrew University of Jerusalem, Israel. The material in this paper was presented in part at the 2016 International Conference on the Science of Electrical Engineering, Eilat, Israel, Nov, 2016, and at the Annual Allerton Conference on Communication, Control, and Computing, Monticello, IL, Sep. 2017. This work was supported in part by the Ministry of Science and Technology, Israel, in conjunction with the Ministry of Europe and Foreign Affairs and with the Ministry of Higher Education, Research, and Innovation, France, in part by the ISF under Grant 1555/18, and in part by the HUJI Federmann Cyber Security Center in conjunction with the Israel National Cyber Directorate (INCD) in the Prime Minister's Office.

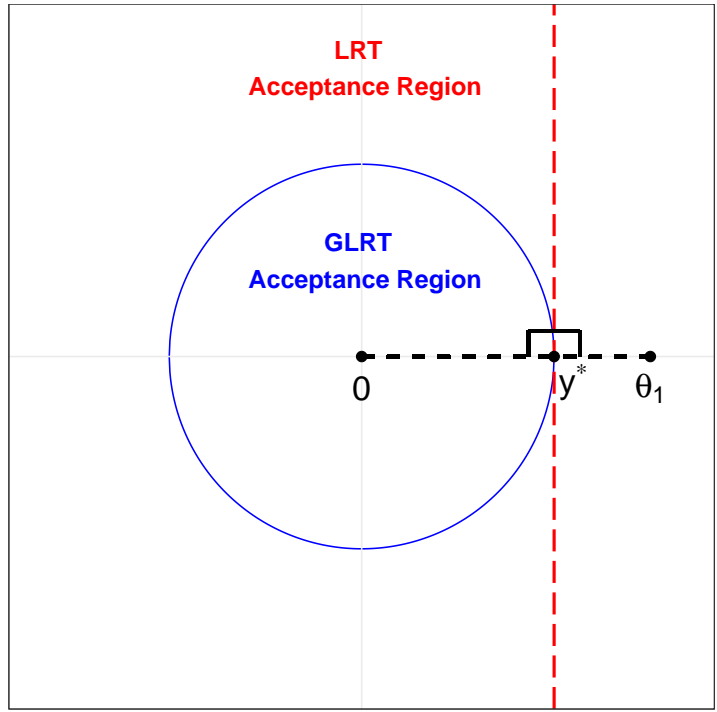

Fig. 1. The acceptance regions for the LRT and GLRT. Both tests achieve the same error exponents which are determined by $y^*$ under both hypotheses.

## I. INTRODUCTION

In the binary hypothesis testing (HT) problem, the goal is to find the optimal tradeoff between the error probabilities given the two hypotheses, and the family of tests which achieves this tradeoff. When both hypotheses are simple, the answer is well known: the Neyman-Pearson lemma states that the optimal test is the (possibly randomized) likelihood-ratio test (LRT), the error-probabilities tradeoff can be computed, and the asymptotics for i.i.d. sequences is given by well-known exponents, see e.g. [1]. However, when at least one of the hypotheses is composite, the situation becomes much more involved. For example, consider the normal location problem (NLP): we have $n$ i.i.d. normal vector measurements $Y_i \sim \mathcal{N}_k(\theta, I)$, where $\theta \in \mathbb{R}^k$. Under hypothesis $\mathcal{H}_0$, $\theta$ is an all-zeros vector, while under $\mathcal{H}_1$ it may be any vector $\theta_1 \in \mathbb{R}^k$ that is non-zero. We will see in the sequel that this example serves as a proxy for richer models.

For known $\theta_1$, the LRT checks the position of the mean $\bar{y}$ with respect to a hyperplane. For unknown $\theta_1$, as we have no reason to prefer any direction with respect to the zero point, a test that "makes sense" should compare $\|\bar{y}\|$ to a threshold. This spherical family of tests happens to coincide with the generalized likelihood-ratio test (GLRT). In Figure1 we depict the LRT and GLRT, with the thresholds chosen such that the hyperplane is tangent to the sphere.

The analysis of this simple model reveals two fundamental problems:

- It is not clear how to define a universal optimality criterion. Specifically, a minimax criterion deems all tests equally useless, as $\theta_1$ may be arbitrarily close to zero.
- Exponential asymptotics are not suitable for quantifying the cost of universality. Specifically, for any LRT there exists a GLRT that yields the same exponents. To see that, consider the regions in Figure1: Under both tests, the same point $y^*$ achieves the minimum distance from zero within the rejection region, as well as the minimum distance from $\theta_1$ within the acceptance region. As for well-behaved regions the error exponents are dictated by the closest point,[1] universality has no cost in terms of exponents.

In the rest of the introduction we outline our framework which solves these two problems, and then describe how we use this framework to provide quantitative measures to the effects of universality, training data, dimensionality and blocklength.

## A. Optimality Criterion

A very large body of works considers the problem of composite hypothesis testing and different notions of universality of a test were suggested, see e.g. [2]. A minimax test considers the worst case among possible distributions. However, this is an overly pessimistic approach, which does not reflect the performance if the distribution happens to be "easier". In the NLP, it assigns a trivial performance to any test. An alternative that overcomes this drawback is uniformly most powerful tests; however, they are only well suited for one-sided tests (in the NLP context, for dimension $k=1$ and where the sign of $\theta_1$ is known).

A partial remedy is given by the notion of *competitive minimax*, introduced in the context of hypothesis testing by Feder and Merhav [3]: for each possible $\theta_1$, we "compete" against the optimal performance when $\theta_1$ is a-priori known, achievable by the LRT; a competitive minimax test maximizes the minimal ratio between the error probability and the LRT error probability, over all possible $\theta_1$ (assuming a uniform prior on the hypotheses). Noticing that for the NLP, the LRT performance is set only by the radius $\|\theta_1\|$, we can see the radius as a measure of "hardness" (a smaller radius means a harder parameter). Competitive minimax indeed means that a test should be "not too bad" for any hardness level. However, two tests that have the same performance for the "critical hardness" will be still considered equivalent, even if for other hardness levels one of them performs better.

---

[1]For discrete alphabets a similar result holds by Sanov's Theorem. In the Gauussian case it can be derived by using "Gaussian types", or by a straightforward approach of bounding the $\mathcal{X}^2$ tail probability.

We present a stronger version, coined *universal minimax* (UMM) optimality, which avoids a Bayesian prior on the hypotheses and ensures optimality at all hardness levels. To that end, we consider the error-probabilities curve under the LRT. That is, let $p_{\text{FA}}$ and $p_{\text{MD}}$ be the error probabilities under $\mathcal{H}_0$ and $\mathcal{H}_1$, respectively. For a family of tests with a varying threshold, the curve is the tradeoff $p_{\text{MD}}(p_{\text{FA}})$ obtained. Comparing the LRT curve for the NLP with two parameters, we notice that the curve for $\theta_1$ with a larger radius dominates that of one with a smaller radius - thus we have a hardness relation that does not require a prior. Now we define *minimax sets* as sets with some bounded hardness (for the NLP, sets of the form $\|\theta_1\| \geq r$), and we require minimax optimality for *all* minimax sets. It is not hard to verify (and we do so in the sequel), that the spherical family of tests (the GLRT) are UMM for the NLP.

Using this definition we can now quantify the "cost of universality", that is, the degradation in performance due to $\mathcal{H}_1$ being composite. As we shall see in Section III, for the NLP with blocklength $n$ the Neyman-Pearson and UMM curves are given by:

$$Q^{-1}(p_{\text{FA}}) + Q^{-1}(p_{\text{MD}}) = d \tag{1a}$$

$$Q^{-1}_{(k),0}(p_{\text{FA}}) = Q^{-1}_{(k),d^2}(1 - p_{\text{MD}}), \tag{1b}$$

respectively, where $Q^{-1}(\cdot)$ and $Q^{-1}_{(\cdot),\cdot}(\cdot)$ are the inverse tail functions of the Gaussian distribution and non-central $\mathcal{X}^2$ distribution, respectively, and $d^2 = n\|\theta_1\|^2$.

Of course, the NLP is of limited interest by itself. We will show that for the wide class of *locally asymptotically normal* (LAN) models, (1) serves as an asymptotic approximation (although the UMM test is not necessarily the GLRT). This class includes, among other examples, testing whether a discrete-alphabet i.i.d. process has some given marginal distribution, and testing whether a Gaussian autoregressive process has some given parameters.

## B. Asymptotics

Consider a sequence of problems with increasing blocklength $n$, some parameter $\Delta$ which controls the "hardness" of the problem at a fixed $n$, and an error probability $p_e$. The most basic approaches to the asymptotic tradeoff between $n$, $\Delta$ and $p_e$ are *multiplicative* and *additive* approximations of $p_e$. Under these two approximations, we seek asymptotic expressions of the forms:

$$p_e = g(\Delta, n)(1 + o(1))$$

$$p_e = g(\Delta, n) + o(1),$$

respectively.[2] As the first is most meaningful for very small $p_e$ and the second for moderate values of $p_e$, it is convenient to consider fixed $\Delta$ and fixed $p_e$, respectively. The multiplicative approximation gives rise to exponential analysis:

$$p_e(n) = e^{-n[E(\Delta) + o(1)]},$$

while the additive approximation gives rise to second-order analysis, which may take the form:

$$\sqrt{\frac{V}{n}} Q^{-1}(p_e) = \Delta(n) + o(1),$$

where the constants $E(\cdot)$ and $V$ are set by the problem parameters. For the latter approximation, we may define $d = \sqrt{n/V}\Delta(n)$ to be the absolute hardness (comparable between different blocklengths) and then fixed error probability means asymptotically fixed hardness, i.e.,

$$Q^{-1}(p_e) = d + o(1).$$

In an information-theoretic context, such analysis is applied to source and channel coding, where $\Delta$ is the gap between the rate and the capacity or rate-distortion function, and $p_e$ is the error probability of the communication or compression scheme. Exponential analysis has been the focus for many year, while second-order analysis gained much popularity in the last decade following [4], with $V$ being the channel or source dispersion.

Let us consider these two approaches in the context of the error-probabilities tradeoff of HT. The multiplicative approximation leads to error exponents. For simple hypotheses, the Stein exponent gives the optimal decay of one of the error probabilities, while the Chernoff exponent gives the optimal balanced decay, see e.g. [1]. However, we noticed above that for the NLP, the same exponents are achievable without knowing $\theta_1$. The same holds for many other parametric (finite-dimensional) models,[3] although, the correction term is affected of course. We suggest thus to turn to the additive approximation, which well-describes the error-probabilities tradeoff for moderate probabilities and also allows to quantify the cost of universality.

---

[2]The first is sometimes relaxed to $\log p_e = g(\Delta, n)(1 + o(1))$.

To that end, consider any LAN model. Using the Fisher information matrix we can transform the problem into an equivalent NLP, and then apply the LRT (for known $\theta_1$) or the UMM (spherical) test (for unknown $\theta_1$) to that model.[4] The resulting performance will be given by (1) up to an additive vanishing correction term, with an appropriate distance $d$, to be specified in the sequel. For example, for a smooth i.i.d LAN model, we have:

$$d = \sqrt{n} \left\| J_{\theta_0}^{1/2}(\theta_1 - \theta_0) \right\|,$$

where $J_\theta$ is the single-measurement Fisher information matrix of the model (In other LAN models, the scaling is not necessarily $\sqrt{n}$). Under an appropriately-defined asymptotic UMM notion, we show that (1b) indeed gives the asymptotic optimal universal performance, thus we can also quantify the asymptotic cost of universality for LAN models.

*C. Further Results*

Using the UMM framework, we can now answer some further questions. First, what happens in the composite problem when training data is added? Specifically, consider that the decoder is given an additional independent sequence of measurements of blocklength $\rho n$, labeled as pertaining to the (unknown) $\theta_1$. We define an extension of the UMM criterion to this case, and evaluate the UMM performance by finding the optimal test. We derive the optimal error-probabilities tradeoff curve under this criterion; indeed, this curve is always between these of (1a) and (1b), and tends to the former or latter when $\rho$ tends to infinity or to zero, respectively. Thus, we also quantify the asymptotic value of training data.

Finally, we explore the tradeoff between blocklength, dimension and training. That is, we allow $k$ and $\rho$ to vary with the blocklength $n$, and consider the interplay between all. We find that the asymptotic UMM performance when $n$ is high, and in addition either the dimension $k$ or training-data ratio $\rho$ are high,[5] we have similar to (1a):

$$Q^{-1}(p_{\text{FA}}) + Q^{-1}(p_{\text{MD}}) = \mathcal{E},$$

where $\mathcal{E}$ is a function of $d$, $k$ and $\rho$ which we explicitly give. We see that in this limit the simple sum of inverse tail functions describes the performance requirements (thus, if one seeks a scalar performance measure, this may be better justified than the common "area under the curve"), while $\mathcal{E}$ describes the "hardness" of the problem. Using this analysis we find, for example, that for the NLP with $n$ i.i.d. measurements, if the dimension grows while $\rho$ is fixed, then the blocklength $n$ must grow as $\sqrt{k}$ in order to keep the UMM performance fixed, both with and without training data.

---

[3]This may not be the case for richer families of distributions, see [5].

[4]Although the UMM test is the GLRT for the NLP, the resulting test is in general not the GLRT for the original LAN model.

*D. Related Work*

We now mention some related work, and highlight the difference from this paper.

Following [3], a non-Bayesian definition that is even closer in spirit to this paper appears in [6]. Also, in a classical work by Gutman [7], the problem of testing $M$ hypotheses, where the source is $K$-th order Markov with training sequences is considered. however, these works remain in the exponential regime.

Refinements of the multiplicative approximation (exponential behavior) can be found in [8], [9]. Such a refinement of Gutman's result is found in [10]. These do improve the estimate of moderate error probabilities, yet they offer very different approximation than our approach which abandons the multiplicative approximation altogether. In [11], the cost of composite testing is considered in a regime similar to ours, but the composite class consists of a finite number of known distributions.

Universal outlier detection was the subject of a recent line of works, see [12] and references therein. However, the approach taken in that line of works is fundamentally different: A set of sequences is jointly classified, and the identity of the outlier sequence(s) is considered a "digital message", in the sense that any misclassification of a sequence constitutes an error event. We, to the contrary, are interested in per-sequence error probabilities.

---

[5]Some further conditions on the relative rate of growth of the parameters are needed, which depend on the LAN model at hand, as we discuss in the sequel.

In [13], binary hypothesis testing with training data for discrete memoryless sources is studied in a game-theoretic framework. The distributions under the hypotheses are known only through training sequences, where one of the sequences has been modified by an adversary. In [14], prediction with training data is considered, where a minimax criterion is used with respect to distributions close to the training-data distribution.

We finally note a recent work [15] which presents a comprehensive framework for learning and inference using local geometry.

*E. Paper Organization*

The rest of this paper is organized as follows. In Section II we introduce the basic definitions for a single-observation (blocklength 1) version of the problem and in Section III we show the optimal performance for the NLP under these definitions with and without training data. In Section IV we add blocklength to the picture and define asymptotic optimality criteria. Section V contains the main result of this paper, showing the asymptotically optimal test for LAN models. Finally in Section VI we analyze sequences of problems with growing dimension.

## II. DEFINITIONS AND PRELIMINARIES

*A. Basic Notation*

Throughout, we use capital letters for random variables, small letters for instances and calligraphic for corresponding alphabets, e.g., $y \in \mathcal{Y}$ is an instance of $Y$. For simplicity, distributions are always either discrete or (vector) continuous. For an event $\mathcal{B}$, the indicator function of the event will be denoted by $1_{\{\mathcal{B}\}}$, i.e. $1_{\{\mathcal{B}\}} = 1$ if $\mathcal{B}$ occurs and $1_{\{\mathcal{B}\}} = 0$ otherwise. The complement of $\mathcal{B}$ will be denoted by $\mathcal{B}^c$. For a random vector $W$ the conditional distribution of $W$ given $\mathcal{B}$ will be denoted by

$$P(W \mid \mathcal{B}).$$

We use $\sim$ for "distributed as", e.g., $X \sim L$, where $L$ is some distribution, discrete or continuous. This may be used also with conditioning, e.g., $X|\mathcal{B} \sim L$ or

$$W \mid X \sim L.$$

If $L$ is a well-known distribution we use standard notation. Most commonly used is the $k$-dimensional normal distribution with mean $\mu$ and covariance matrix $\Sigma$, denoted by $\mathcal{N}_k(\mu, \Sigma)$. For $k=1$ we will omit the subscript $k$, e.g. $\mathcal{N}(0,1)$.

If a sequence of random vectors $W_n$ converges in distribution to $W$, it will be denoted by $W_n \rightsquigarrow W$. If $W$ has a distribution $L$, or a distribution with a standard notation, such as $\mathcal{N}(0,1)$, then we may also denote it by $W_n \rightsquigarrow L$ or $W_n \rightsquigarrow \mathcal{N}(0,1)$.

The notation $o_p(1)$ is short for a sequence of random vectors that converges to zero in probability. The expression $O_p(1)$ denotes a sequence that is bounded in probability as defined next.

*Definition 1:* A sequence of random vectors $W_n$ is bounded in probability if for every $\varepsilon > 0$ there exists $M$ such that
$$\sup_n P(\|W_n\| > M) < \varepsilon.$$
We will denote this by $W_n = O_p(1)$.

For deterministic vectors we write $o(1)$ or $O(1)$.

We will use a subscript to indicate distributions indexed by a parameter. That is, $(\Pi_\theta | \theta \in \Theta)$ is a family of distributions (or model) indexed by a parameter set $\Theta \subseteq \mathbb{R}^k$. For brevity of notation when the context is clear we simply write $\Pi_\theta$. For an event $\mathcal{B}$ we denote by
$$P(\mathcal{B}; \theta)$$
the probability of the event if $Y \sim \Pi_\theta$. Expectation is similarly denoted.

## B. Discriminant Rules and Error Probabilities

Let $(\Pi_\theta | \theta \in \Theta)$ be a family of distributions on some alphabet $\mathcal{Y}$ parametrized by $\theta \in \Theta \subseteq \mathbb{R}^k$. Suppose that $Y \sim \Pi_\theta$ and consider the following problem, of testing between known $\theta_0$ and unknown $\theta_1$:[6]

$$\begin{aligned} \mathcal{H}_0 &: Y \sim \Pi_{\theta_0} \\ \mathcal{H}_1 &: Y \sim \Pi_{\theta_1} \text{ for some } \theta_1 \in \tilde{\Theta}, \end{aligned} \quad (2)$$

---

[6] Unlike the introduction where we considered a sequence of $n$ measurements, here we consider a single one. We will re-introduce blocklength in Section IV.

where $\tilde{\Theta} = \Theta \setminus \{\theta_0\}$ is the set of possible parameters under $\mathcal{H}_1$ which we refer to as the set of alternative parameters or alternatives.

A discriminant rule (detector) without training sequence is a (possibly randomized) statistic, $R(y)$, taking values in $\{0, 1\}$. Given the element $\theta_1$, it induces the "false-alarm" and "missed-detection" error probabilities:

$$p_{\text{FA}}(R) = P(R(Y) = 1; \theta_0),$$
$$p_{\text{MD}}(R; \theta_1) = P(R(Y) = 0; \theta_1). \tag{3}$$

We will also refer to (3) as the missed-detection probability against $\theta_1$. In statiatical literature the false-alarm and missed-detection are often referred to as type-I and type-II errors respectively. We denote by

$$\mathcal{R}(p) = \{R \mid p_{\text{FA}}(R) \leq p\}$$

the set of all detectors with level of significance $p$. Note that the level of significance is guaranteed regardless of $\theta_1$. The optimal performance for a given $\theta_1$ is given by:

$$p^*_{\text{MD}}(p_{\text{FA}}; \theta_1) = \min_{R \in \mathcal{R}(p_{\text{FA}})} p_{\text{MD}}(R; \theta_1).$$

By the Neyman-Pearson lemma, this optimum is given by the likelihood ratio test (LRT): $R(y) = 0$ when

$$T < \log \frac{L_{\theta_0}(y)}{L_{\theta_1}(y)}, \tag{4}$$

where $L_\theta(\cdot)$ is the likelihood function and $T \in \mathbb{R}$ is a threshold; for discrete distributions the likelihood is simply the probability mass function of $\Pi_\theta$ evaluated at $y$ and for continuous distributions it is the probability density function evaluated at $y$. Scanning over the threshold $T$ gives the whole tradeoff curve. In general, randomization between LRT with different thresholds is needed, but this will be ignored in the sequel as it is insignificant in the regimes of interest.

In our compound setting (2), where $\theta_1$ is unknown, a commonly used detector is the generalized likelihood-ratio test (GLRT): $R(y) = 0$ when

$$T < \log \frac{L_{\theta_0}(y)}{L^*_{\tilde{\Theta}}(y)},$$

where $L^*_{\tilde{\Theta}}(y) = \sup_{\theta_1 \in \tilde{\Theta}} L_{\theta_1}(y)$.

## C. Universal Minimax Optimality

In order to define our universal optimality criterion for the problem (2) we require a "hardness" ordering of *all* alternatives with respect to the LRT performance, $p^*_{\text{MD}}$, which is independent of $p_{\text{FA}}$. This ordering is defined as follows.

*Definition 2:* Let $\theta, \theta' \in \tilde{\Theta}$, $\theta \neq \theta'$. We say that $\theta'$ is harder than $\theta$ if for all $0 < p_{\text{FA}} < 1$

$$p^*_{\text{MD}}(p_{\text{FA}}; \theta) < p^*_{\text{MD}}(p_{\text{FA}}; \theta'). \tag{5}$$

We denote this relation by $\theta' \prec \theta$. If the inequality in (5) can be replaced by equality or weak inequality, we denote the relation by $\theta' \simeq \theta$ ("equally hard") or $\theta' \preccurlyeq \theta$ ("at least as hard"), respectively. Note that this is not necessarily a partial order relation as it may not be anti-symmetric.

*Definition 3:* If for any two alternatives $\theta, \theta' \in \tilde{\Theta}$ either $\theta' \preccurlyeq \theta$ or $\theta \preccurlyeq \theta'$, then the model $\Pi_\theta$ is said to be degraded at $\theta_0$.

We define sets of alternative for which the hardness of detection (when $\theta_1$ is known) is below some threshold:

*Definition 4:* For an alternative $\theta \in \tilde{\Theta}$ denote by $\tilde{\Theta}_\theta$ the set of all alternatives which are not harder than $\theta$:

$$\tilde{\Theta}_\theta = \left\{ \theta_1 \in \tilde{\Theta} \middle| \theta \preccurlyeq \theta_1 \right\}.$$

We call such sets minimax sets.

Finally we define a strong optimality criterion, which requires minimax optimality within *any* minimax set.

*Definition 5:* Suppose that $\Pi_\theta$ is degraded at $\theta_0$. $R^*$ is said to be a universal minimax (UMM) detector with significance level $p_{\text{FA}} \in (0,1)$, if $R^*$ minimizes $\sup_{\theta_1 \in \tilde{\Theta}_\theta} p_{\text{MD}}(R; \theta_1)$ among all detectors in $\mathcal{R}(p_{\text{FA}})$, simultaneously for all $\theta \in \tilde{\Theta}$.

In order to elucidate on the meaning of UMM optimality, we revisit the comparison with competitive minimax, introduced in Section I-A, in the context of degraded families. Consider all detectors that have some fixed false-alarm probability $p_{\text{FA}}$ for all $\theta$. By the definition of degradedness, we may characterize the hardness of a parameter $\theta$ by the corresponding LRT performance $p^*_{\text{MD}}(p_{\text{FA}}; \theta)$. For an arbitrary rule and for any $0 < p^* < 1$, let

$$p_{\text{MD}}(p^*) = \max_{\theta: p^*_{\text{MD}}(p_{\text{FA}}; \theta) = p^*} p_{\text{MD}}(p_{\text{FA}}; \theta).$$

A competitive minimax rule, is a rule which minimizes the worst-case ratio:

$$\max_{p^*} \frac{p_{\text{MD}}(p^*)}{p^*}.$$

There is no requirement for optimality of $p_{\text{MD}}(p^*)$ for values of $p^*$ other than the worst-case one. UMM optimality, on the other hand, typically means that the rule minimizes $p_{\text{MD}}(p^*)$ simultaneously for all $p^*$. Indeed this is a very strict requirement, but we will see that it can be satisfied in the simple normal location problem, and that in an asymptotic sense to be defined later, it can be satisfied for many other problems as well.

## D. Discriminant Rules with a Training Sequence

Let $X$ be independent of $Y$ and labeled as coming from a "related" distribution with the same alternative parameter $\theta_1$. Specifically, $X \sim \tilde{\Pi}_{\theta_1}$, where $\tilde{\Pi}_\theta$ is a family of distributions on some alphabet $\mathcal{X}$. The distributions $\Pi_\theta$ are defined for the same set $\tilde{\Theta}$ as above, but they may not be identical to $\tilde{\Pi}_\theta$.[7] Given $\theta_1$, $X$ is independent of $Y$. We refer to $X$ as a "training sequence" or "training data" for the alternative. Our main interest will be in discriminant rules for this problem which use "training data". In section IV we will consider sequences of such problems.

For events $\mathcal{B}$ that depend upon both $Y$ and $X$, we use

$$P(\mathcal{B}; \theta_y, \theta_x)$$

to denote the probability when $Y \sim \Pi_{\theta_y}$ and $X \sim \tilde{\Pi}_{\theta_x}$. Similarly, for the conditional probability of $\mathcal{B}$ given $X$ we use

$$P(\mathcal{B} \mid X; \theta_y, \theta_x).$$

If only one parameter is denoted then it is $\theta_y$.

When training data is added to the problem 2, it can be thought of as an HT problem over the alphabet $\mathcal{X} \times \mathcal{Y}$ of the following form:

$$\begin{aligned} \mathcal{H}_0 &: (X,Y) \sim \tilde{\Pi}_{\theta_1} \times \Pi_{\theta_0} \\ \mathcal{H}_1 &: (X,Y) \sim \tilde{\Pi}_{\theta_1} \times \Pi_{\theta_1}, \end{aligned} \tag{6}$$

---

[7]For example, $X$ and $Y$ may be i.i.d. with the same marginal but with a different blocklength, such that for some base alphabet $\mathcal{A}$, $\mathcal{Y} = \mathcal{A}^n$ but $\mathcal{X} = \mathcal{A}^m$.

for some $\theta_1 \in \tilde{\Theta}$. A discriminant rule (detector) with a training sequence is a (possibly randomized) statistic $R_x(y)$ of $(x, y)$ taking values in $\{0, 1\}$. It induces *conditional* error probabilities

$$p_{\text{FA}}(R_x|x) = P(R_x(Y) = 1 \mid X = x; \theta_0)$$

$$p_{\text{MD}}(R_x|x; \theta_1) = P(R_x(Y) = 0 \mid X = x; \theta_1).$$

The (average) missed-detection probability is denoted by

$$p_{\text{MD}}(R_x; \theta_1) = E\left[p_{\text{MD}}(R_x|X; \theta_1)\right]$$

$$= P(R_X(Y) = 0; \theta_1, \theta_1).$$

Now, we denote as

$$\underline{\mathcal{R}}(p) = \{R_x \mid p_{\text{FA}}(R_x|x) \leq p, \forall x \in \mathcal{X}\}$$

the set of detectors with guaranteed significance level $p$ for every $x$. The corresponding UMM criterion is as follows.

*Definition 6:* Suppose that $\Pi_\theta$ is degraded at $\theta_0$. $R_x^*$ is said to be a universal minimax detector with a training sequence and guaranteed significance level $p_{\text{FA}} \in (0, 1)$, if it minimzes $\sup_{\theta_1 \in \tilde{\Theta}_\theta} p_{\text{MD}}(R_x; \theta_1)$ among all detectors in $\underline{\mathcal{R}}(p_{\text{FA}})$, simultaneously for all $\theta \in \tilde{\Theta}$.

Note that, while the missed-detection probability is averaged over the random $X$, with guaranteed significance level we require that the false-alarm probability will be at most $p_{\text{FA}}$ for *all* $x$. In Section VI-D we will discuss UMM with average significance level.

## III. UNIVERSAL MINIMAX FOR THE NORMAL LOCATION PROBLEM

We now address the simple NLP, which will be used in the sequel as an asymptotic proxy for more involved problems. Notice that for now we do not introduce blocklength, i.e., there is only a single measurement. Let the parameter set be $\Theta = \mathbb{R}^k$. For the NLP we shall denote the parameters by $\mu$ instead of $\theta$ to distinguish this problem from the problems in later sections. Consider the problem (2) with $\Pi_\mu = \mathcal{N}_k(\mu, \Sigma)$, for some known positive definite covariance matrix $\Sigma$. Without loss of generality we assume that $\mu_0 = 0$ and $\Sigma = I$, otherwise we can standardize $Y$. To be more specific, we choose some diagonalization $\Sigma = A\Lambda A^t$ for an orthogonal $A$, where $\Lambda = \text{diag}(\lambda_1, \ldots, \lambda_k)$, and we define the corresponding square-root matrix

$$\Sigma^{1/2} = A\Lambda^{1/2}, \tag{7}$$

where $\Lambda^{1/2} = \text{diag}\left(\sqrt{\lambda_1}, \ldots, \sqrt{\lambda_k}\right)$. We can standardize $Y$ by replacing it with $\Sigma^{-1/2}(Y - \mu_0)$ which is distributed $\mathcal{N}_k(0, I)$, under $\mathcal{H}_0$. We denote the Mahalanobis distance between the null and alternative distribution by:[8]

$$\Delta^2 = \|\mu_1\|^2.$$

If we knew $\mu_1$ we would be able to use the LRT (4); although the resulting test is well known, we include a short derivation. The log likelihood-ratio of the normal vector $Y$ at $y$ is given by

$$\log \frac{L(y; \mu_1, I)}{L(y; 0, I)} = \mu_1^t y - \frac{\Delta^2}{2}, \tag{8}$$

where $L(y; \mu, I)$ denotes the probability density function (likelihood) of $\mathcal{N}_k(\mu, I)$ and $\mu_1^t$ is the transpose of $\mu_1$. The log likelihood-ratio (8) is normally distributed with mean $-\Delta^2/2$ or $\Delta^2/2$ under $\mathcal{H}_0$ or $\mathcal{H}_1$, respectively, and variance $\Delta^2$ under both hypotheses. The LRT (4) amounts to the family of rules with acceptance regions given by

$$\mu_1^t y < T,$$

for $T \in \mathbb{R}$. Using the distribution of (8), it is easy to see that the error probabilities of the LRT are given by:

$$p_{\text{FA}}(R) = Q\left(\frac{T}{\Delta}\right)$$
$$p_{\text{MD}}(R) = 1 - Q\left(\frac{T - \Delta^2}{\Delta}\right),$$

where $Q(\cdot)$ is the tail distribution function of the standard normal distribution. Equating the threshold $T$, we find that $p_{\text{FA}}$ and $p_{\text{MD}}$ satisfy

$$Q^{-1}(p_{\text{FA}}(R)) + Q^{-1}(p_{\text{MD}}(R)) = \Delta. \tag{9}$$

From (9) it is clear that the model $\Pi_\mu$ is degraded at $\mu_0$. Furthermore, since the hardness relation (inversely) agrees with the norm, the minimax sets $\tilde{\Theta}_\mu$ depend only on the norm, and are given by all $\tilde{\Theta}(d)$, $d > 0$, where

$$\tilde{\Theta}(d) = \left\{\mu_1 \in \tilde{\Theta} \,\middle|\, \|\mu_1\| \geq d\right\}. \tag{10}$$

---

[8]Notice that it is also twice the KL-divergence between $\Pi_{\mu_1}$ and $\Pi_0$.

We notice that sets of this form seem like a natural way to keep $\mu_1$ separated from $\mu_0$ and thus facilitate non-trivial minimax analysis. Indeed, many works have considered such sets. For example, in [16], $\mu_1$ is restricted to satisfy $\|\mu_1 - \mu_0\| \geq r$, for some $r > 0$ and some metric $\|\cdot\|$; taking the metric to be the Euclidean distance, recovers our minimax sets. The UMM framework gives an operational meaning to this classical analysis in terms of LRT "hardness", and allows to generalize it to further problems.

## A. Universal minimax without training data

For unknown $\mu_1$, the commonly used GLRT amounts to the family of rules with acceptance regions in the form of spheres around the origin: $\|y\|^2 < T$, where for a false-alarm probability $p_{\text{FA}}$ we set

$$T = Q_{(k)}^{-1}(p_{\text{FA}}), \qquad (11)$$

where $Q_{(k)}^{-1}(\cdot)$ is the inverse tail function of the $\mathcal{X}^2$ distribution with $k$ degrees of freedom. The performance of the GLRT is given by $p_{\text{MD}}$ satisfying

$$Q_{(k)}^{-1}(p_{\text{FA}}) = Q_{(k),\Delta^2}^{-1}(1 - p_{\text{MD}}), \qquad (12)$$

where $Q_{(k),\lambda}^{-1}(\cdot)$ is the inverse tail function of the non-central $\mathcal{X}^2$ distribution with $k$ degrees of freedom and non-centrality parameter $\lambda$.

*Proposition 1:* The GLRT (11) is a universal minimax detector with significance level $p_{\text{FA}}$.

*Proof:* The GLRT with any threshold $T$ is minimax over the minimax set $\tilde{\Theta}(\Delta)$, for all $\Delta$ and with minimax performance given by (12), see [2, Problem 8.29]. Since $\|y\|^2$ has a $\mathcal{X}^2$ distribution with $k$-degrees of freedom under $\mathcal{H}_0$, then (11) ensures the significance level $p_{\text{FA}}$. ∎

## B. Universal minimax with training data

Suppose that $\tilde{\Pi}_\mu = \mathcal{N}_k\left(\mu, \frac{1}{\rho}\Sigma\right)$, for some positive known scalar $\rho$ which can be thought of as the "quality" of the training data.[9] As without training data, without loss of generality we take $\mu_0 = 0$ and $\Sigma = I$, otherwise we standardize both $X$ and $Y$. Specifically, we seek a UMM rule for the problem

$$\begin{aligned} \mathcal{H}_0 &: (X, Y) \sim \mathcal{N}_k\left(\mu_1, \frac{1}{\rho}I\right) \times \mathcal{N}_k(0, I) \\ \mathcal{H}_1 &: (X, Y) \sim \mathcal{N}_k\left(\mu_1, \frac{1}{\rho}I\right) \times \mathcal{N}_k(\mu_1, I). \end{aligned} \quad (13)$$

---

[9]This choice may seem arbitrary at the moment. However, consider multiple measurements, then $X$ and $Y$ can stand for sample averages pertaining to the same marginal distribution with different blocklength.

As a means for solution we consider a Bayesian variant of (13) where $\mu_1$ itself has a uniform prior on a sphere of radius $\Delta$ centered at the origin, see, e.g., [17, Chapter 15, problem 4] for a similar approach.[10] To be more specific, let $M$ be uniformly distributed on the $\Delta$-sphere in $\mathbb{R}^k$ and suppose that the conditional distribution is

$$X \mid \{M = m\} \sim \mathcal{N}_k\left(m, \frac{1}{\rho}I\right).$$

Notice that this distribution is parametrized by a single parameter, $\Delta$. A similar assumption holds for $Y$, and we can state an HT problem similar to (13) as follows:

$$\begin{aligned} \mathcal{H}_0 &: (X, Y) \mid \{M = m\} \sim \mathcal{N}_k\left(m, \frac{1}{\rho}I\right) \times \mathcal{N}_k(0, I) \\ \mathcal{H}_1 &: (X, Y) \mid \{M = m\} \sim \mathcal{N}_k\left(m, \frac{1}{\rho}I\right) \times \mathcal{N}_k(m, I), \end{aligned} \quad (14)$$

where under both hypotheses $M$ is uniform over the $\Delta$-sphere. The next lemma, which is proven in appendix A, gives the LRT for this problem.

*Lemma 1:* For any $\Delta > 0$ the acceptance region of the LRT for the problem (14) is given by

$$c_k(\Delta \|\rho x + y\|) > T \cdot c_k(\Delta \|\rho x\|), \quad (15)$$

where $T > 0$ and

$$c_k(\tau) = \frac{\tau^{(\frac{k}{2}-1)}}{(2\pi)^{\frac{k}{2}} \cdot \mathcal{I}_{(\frac{k}{2}-1)}(\tau)} \quad (16)$$

is the normalizing constant of the Von Mises-Fisher distribution, with $\mathcal{I}_q(\cdot)$ denoting the modified Bessel function of the first kind of order $q$.

The Von Mises-Fisher distribution arises from the principle of maximum entropy which is known to be closely related to minimax optimality, see e.g., [14]. To be more specific, when $M$ has a uniform prior on the $\Delta$-sphere, the posterior distribution of $\frac{1}{\Delta}M \mid X = x$ is the Von Mises-Fisher distribution with mean direction $x/\|x\|$ and concentration parameter $(\rho\Delta\|x\|)$ [18], which is the maximum entropy distribution on the unit-sphere with a (non-zero) constraint on the first moment, see e.g., [19, Problem 15.3.9].

---

[10]In fact the GLRT (11) can also be derived as the LRT for this problem proving it is minimax on the set $\tilde{\Theta}(\Delta)$.

Returning to our original problem, the acceptance region in (15) yields a UMM rule with guaranteed significance level as stated in the next theorem which is proven in appendix A.

*Theorem 1:* A universal minimax discriminant rule with a training sequence and guaranteed significance level $p_{\text{FA}}$ for the problem (13) is given by the following acceptance region:

$$\|\rho x + y\|^2 < Q_{(k),\vartheta_0}^{-1}(p_{\text{FA}}), \tag{17}$$

where $\vartheta_0 = \|\rho x\|^2$. We denote this rule by $\underline{R}_{\boldsymbol{x}}$. The missed-detection probability against $\mu_1$ is given by

$$\underline{p}_{\text{MD}}^{\text{umm}}(p_{\text{FA}}, \Delta, \rho, k) = 1 - E\left[Q_{(k),\vartheta_1}\left(Q_{(k),\vartheta_0}^{-1}(p_{\text{FA}})\right); \mu_1\right], \tag{18}$$

where $\vartheta_1 = \|\rho x + \mu_1\|^2$ and $\Delta = \|\mu_1\|$.

Figure 2 demonstrates the acceptance regions of the UMM test (17) for fixed $p_{\text{FA}}$ as $\rho$ grows. As the LHS of (18) suggests, the UMM performance depends upon the parameters only through $\Delta$ and $\rho$. Note that for $\rho = 0$, $\underline{p}_{\text{MD}}^{\text{umm}}(p_{\text{FA}}, \Delta, \rho, k)$ reduces to the GLRT performance in (12), which does not use training data. Intuitively, as $\rho$ tends to $0$ the variance of the training data tends to $\infty$ which makes it useless and thus can be ignored. In Section V we shall see a similar asymptotic relation for a different family of models. Note also that $\underline{p}_{\text{MD}}^{\text{umm}}(p_{\text{FA}}, \cdot, \rho, k)$ is monotonically decreasing and

$$\lim_{\Delta \to \infty} \underline{p}_{\text{MD}}^{\text{umm}}(p_{\text{FA}}, \Delta, \rho, k) = 0.$$

For notational convenience we extend $\underline{p}_{\text{MD}}^{\text{umm}}(p_{\text{FA}}, \cdot, \rho, k)$ to $[0, \infty]$ accordingly. That is, for $\Delta = \infty$ we define $\underline{p}_{\text{MD}}^{\text{umm}}(p_{\text{FA}}, \Delta, \rho, k) = 0$.

We notice that our ability to derive UMM performance hinges on requiring *guaranteed* significance level. If $p_{\text{FA}}$ is allowed to depend on $X$, one may use the "empirical quality" information $\|X\|/\Delta$, thus the optimal tests for different minimax sets (different $\Delta$) may not coincide, see Section VI-D.

## C. Performance Tradeoff Curves

Summarizing the results of this section, for the NLP with any given $\|\mu_1\| = \Delta$ we have the following error-probabilities tradeoffs: the optimal (known $\mu_1$) tradeoff (9), the UMM tradeoff

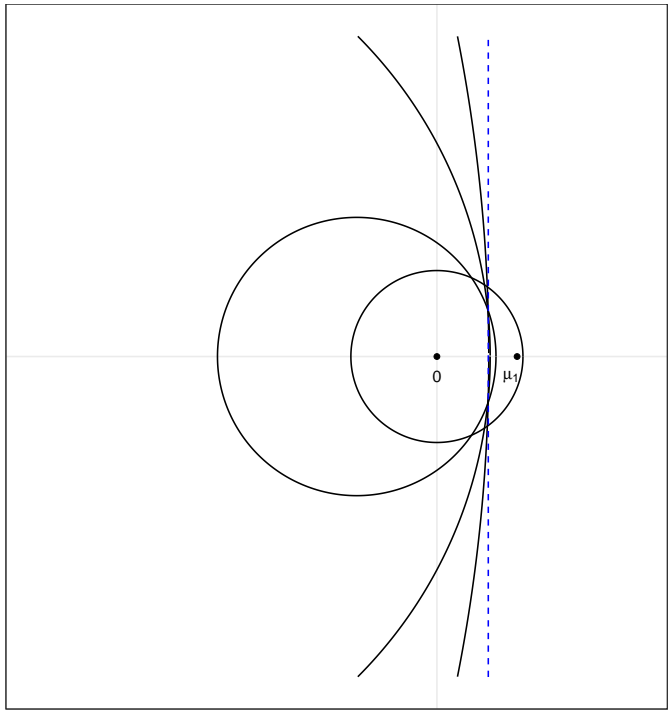

Fig. 2. The UMM acceptance regions, for $k = 2$, $\Delta = 2$. We choose $\mu_1 = (\Delta, 0)$ and $\rho = 0, 1, 5, 20$. $\rho = 0$ corresponds to the GLRT. The radii of the spheres guarantee $p_{\text{FA}} = 0.1$. For simplicity, we only plot the regions for $x = E[x] = \mu_1$. The hyperplane in dashed-blue is the boundary of the acceptance region of the LRT with the same false-alarm probability. As $\rho$ increases the sphere centers migrate left, and at the same time their radii grow, such that locally they approach the LRT boundary.

(12) and the UMM tradeoff with training data and guaranteed significance level (18). As expected, the latter tradeoff is always between the first two, and tends to the first or the second when $\rho$ approaches infinity or zero, respectively. In figure 3 the tradeoff curves are shown to approach the tradeoff curve of the LRT as $\rho$ grows.

## IV. Beyond the NLP: Asymptotic Universal Minimax

Our UMM criteria are very stringent. Indeed, it is difficult to think of problems beyond the NLP where they strictly apply. In this section we define *asymptotic* variants that may be relevant to models of interest, including LAN models, which we introduce in Section V. To that end, we replace the single measurements $Y$ and $X$ by sequences.

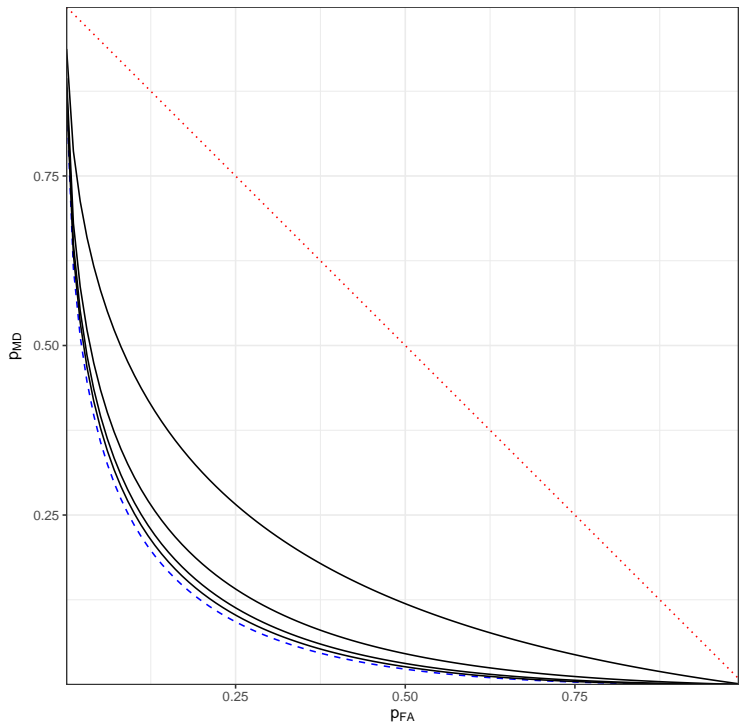

Fig. 3. The parameters $k, \Delta$ and $\rho$ are as in figure 2. The trivial performance (dotted red) and the LRT (dashed blue). Between them are the UMM performances. As $\rho$ grows the UMM performance approaches the LRT.

## A. Sequences of HT problems

In this section we define sequences of HT problems with the corresponding error probabilities. These are straightforward extensions of the problem (2). To be more specific, let $\mathcal{A}$ be some alphabet and denote by $\mathcal{A}^\infty$ the set of all sequences of elements of $\mathcal{A}$ indexed by the natural numbers $\mathbb{N}$. $\mathcal{A}^\infty$ has the same role as $\mathcal{Y}$ in (2). Now $\Pi_\theta$ will be used to denote a family of random processes on $\mathcal{A}^\infty$ with time domain $\mathbb{N}$. We consider test sequences $\mathbf{Y} = (Y_1, \ldots, Y_n)$ which consist of the *first* $n$ samples of the process, as opposed to any successive $n$ samples (this is significant as the models need not be stationary). The distribution of the first $n$ samples from the process $\Pi_\theta$ is denoted by $\Pi_\theta^{(n)}$. Fix some known $\theta_0$ and consider the following problem

$$\begin{aligned}\mathcal{H}_0 &: \mathbf{Y} \sim \Pi_{\theta_0}^{(n)} \\ \mathcal{H}_1 &: \mathbf{Y} \sim \Pi_{\theta_1}^{(n)} \text{ for some } \theta_1 \in \tilde{\Theta}.\end{aligned} \qquad (19)$$

For any $n$ this problem is the same as (2) on the alphabet $\mathcal{A}^n$ with $\Pi_\theta$ replaced by $\Pi_\theta^{(n)}$. Similarly, we have the corresponding problem with training data

$$\begin{aligned} \mathcal{H}_0 &: (\mathbf{X}, \mathbf{Y}) \sim \Pi_{\theta_1}^{(n_X)} \times \Pi_{\theta_0}^{(n)} \\ \mathcal{H}_1 &: (\mathbf{X}, \mathbf{Y}) \sim \Pi_{\theta_1}^{(n_X)} \times \Pi_{\theta_1}^{(n)}, \end{aligned} \quad (20)$$

for some $\theta_1 \in \tilde{\Theta}$ and $n_X \in \mathbb{N}$. For any pair $(n, n_X)$ this problem is the same problem as in (6) on the alphabet $\mathcal{A}^{n_X} \times \mathcal{A}^n$ with $\tilde{\Pi}_{\theta_1}$ replaced by $\Pi_{\theta_1}^{(n_x)}$. Thus if we fix $n$ and $n_X$, all the definitions in section II hold.

*Example 1:* Let $\Pi_\theta$ be the process of i.i.d. $k$-dimensional normal variables with unit variance and let $\theta_0 = 0$. We have $\mathbf{Y} \sim \mathcal{N}_k^n(\theta, I)$ and $\mathbf{X} \sim \mathcal{N}_k^{n_X}(\theta, I)$. Using the same arguments as in section III it can be easily shown that the UMM rule without training data is still the GLRT (11) with $y$ replaced by $\sqrt{n}\bar{y}$, where $\bar{y}$ is the sample mean. The UMM performance is given by (12) with $\Delta = \sqrt{n}\|\theta_1\|$, or $\underline{p}_{\text{MD}}^{\text{umm}}(p_{\text{FA}}, \sqrt{n}\|\theta_1\|, 0, k)$. Similarly, the UMM rule with training data and guaranteed significance level $p_{\text{FA}}$ is given by (17) with $x$ replaced by $\sqrt{n}\bar{x}$ and with $\rho = n_X/n$. The UMM performance is given by $\underline{p}_{\text{MD}}^{\text{umm}}(p_{\text{FA}}, \sqrt{n}\|\theta_1\|, \rho, k)$. In other words, the UMM performance (with and without training data) is the same as the UMM performances for a single observation in III with $\Delta = \sqrt{n}\|\theta_1\|$. Furthermore, notice that unless the sequence $d_n = \sqrt{n}\|\theta_n\|$ is bounded away from $0$ and $\infty$ the UMM performance against $\theta_n$ becomes trivial, i.e. $p_{\text{MD}} = 1 - p_{\text{FA}}$ or $p_{\text{MD}} = p_{\text{FA}} = 0$, respectively.

We will consider *sequences* of the problems in (19) and in (20) indexed by $n$. The alternative $\theta_1$ replaced by a sequence $\theta_n$ and the training-data blocklength is given by a sequence $n_X(n)$. Accordingly, we have a sequence of discriminant rules and corresponding error probabilities. Specifically, if $R^{(n)}$ is a sequence of discriminant rules for (19) then $R^{(n)}$ is a rule for the $n$-th problem and the sequence of error probabilities will be denoted by

$$p_{\text{FA}}^{(n)}\left(R^{(n)}\right) = P\left(R^{(n)}(\mathbf{Y}) = 1; \theta_0\right),$$
$$p_{\text{MD}}^{(n)}\left(R^{(n)}; \theta_n\right) = P\left(R^{(n)}(\mathbf{Y}) = 0; \theta_n\right).$$

Similarly, if $R_{\boldsymbol{x}}^{(n)}$ is a sequence of discriminant rules for (20) the sequence of *conditional* error probabilities will be denoted by

$$p_{\text{FA}}^{(n)}\left(R_{\boldsymbol{x}}^{(n)} \middle| \boldsymbol{x}\right) = P\left(R_{\boldsymbol{x}}^{(n)}(\mathbf{Y}) = 1 \middle| \mathbf{X} = \boldsymbol{x}; \theta_0\right),$$

$$p_{\text{MD}}^{(n)}\left(R_{\boldsymbol{x}}^{(n)}\Big|\boldsymbol{x};\theta_n\right) = P\left(R_{\boldsymbol{x}}^{(n)}(\mathbf{Y}) = 0\Big|\mathbf{X} = \boldsymbol{x};\theta_n\right).$$

The average missed-detection probability is denoted by:

$$\begin{aligned}p_{\text{MD}}^{(n)}\left(R_{\boldsymbol{x}}^{(n)};\theta_n\right) &= E\left[p_{\text{MD}}\left(R_{\boldsymbol{x}}^{(n)}\Big|X;\theta_n\right)\right] \\ &= P\left(R_X^{(n)}(Y) = 0;\theta_n,\theta_n\right).\end{aligned}$$

## B. Asymptotic UMM

Based on the definitions in Section IV-A we define the asymptotic optimality criterion and related concepts.

*Definition 7:* Let $\{\theta_n\}$ and $\{\theta'_n\}$ be two sequences of alternatives. We say that the sequence $\{\theta'_n\}$ is harder than $\{\theta_n\}$ if

$$\limsup_{n\to\infty} p_{\text{MD}}^*(p_{\text{FA}};\theta_n) < \limsup_{n\to\infty} p_{\text{MD}}^*(p_{\text{FA}};\theta'_n),$$

for all $0 < p_{\text{FA}} < 1$. We denote this relation by $\{\theta'_n\} \prec \{\theta_n\}$. We define $\{\theta'_n\} \preccurlyeq \{\theta_n\}$ as in definition 2.

*Definition 8:* If for any two sequences of alternatives $\{\theta_n\}, \{\theta'_n\}$ either $\{\theta'_n\} \preccurlyeq \{\theta_n\}$ or $\{\theta_n\} \preccurlyeq \{\theta'_n\}$, then the model $\Pi_\theta$ is said to be asymptotically degraded at $\theta_0$.

For a sequence of alternatives $\{\theta_n\}$ we define asymptotic minimax sets, parallel to definition 4, as follows.

*Definition 9:* Asymptotic minimax sets are given by

$$\tilde{\Theta}_{\{\theta_n\}} = \left\{\{\theta'_n\} \subseteq \tilde{\Theta}\Big|\{\theta_n\} \preccurlyeq \{\theta'_n\}\right\}.$$

The asymptotically UMM optimality criterion is defined as follows.

*Definition 10:* Suppose that the model $\Pi_\theta$ is asymptotically degraded at $\theta_0$ and let $R^{(n)}$ be a sequence of discriminant rules. $R^{(n)}$ is said to be an asymptotically universal minimax (AUMM) rule with significance level $p_{\text{FA}}$, if $\limsup p_{\text{FA}}^{(n)}\left(R^{(n)}\right) \leq p_{\text{FA}}$ and for any other sequence of rules $\tilde{R}^{(n)}$ such that $\limsup p_{\text{FA}}^{(n)}\left(\tilde{R}^{(n)}\right) \leq p_{\text{FA}}$ the following holds

$$\sup_{\{\theta'_n\}\in\tilde{\Theta}_{\{\theta_n\}}}\left[\limsup_{n\to\infty} p_{\text{MD}}^{(n)}\left(R^{(n)};\theta'_n\right)\right] \leq \sup_{\{\theta'_n\}\in\tilde{\Theta}_{\{\theta_n\}}}\left[\limsup_{n\to\infty} p_{\text{MD}}^{(n)}\left(\tilde{R}^{(n)};\theta'_n\right)\right],$$

simultaneously for all $\tilde{\Theta}_{\{\theta_n\}}$.

Finally, we extend definition 10 to discriminant rules with training data as follows.

*Definition 11:* Let $R_{\boldsymbol{x}}^{(n)}$ be a sequence of discriminant rules for (20). $R_{\boldsymbol{x}}^{(n)}$ is said to have a guaranteed asymptotic significance level $p_{\text{FA}}$ if almost surely for every $\boldsymbol{x}$

$$\limsup_{n \to \infty} p_{\text{FA}}^{(n)} \left( R_{\boldsymbol{x}}^{(n)} \Big| \boldsymbol{x} \right) \leq p_{\text{FA}}.$$

*Definition 12:* Suppose that the model $\Pi_\theta$ is asymptotically degraded at $\theta_0$ and let $R_{\boldsymbol{x}}^{(n)}$ be a sequence of discriminant rules for (20). $R_{\boldsymbol{x}}^{(n)}$ is said to be an AUMM with guaranteed significance level $p_{\text{FA}}$ if it has a guaranteed asymptotic significance level $p_{\text{FA}}$ and for any other sequence of rules $\tilde{R}_{\boldsymbol{x}}^{(n)}$ which has a guaranteed asymptotic significance level $p_{\text{FA}}$ the following holds

$$\sup_{\{\theta_n'\} \in \tilde{\Theta}_{\{\theta_n\}}} \left[ \limsup_{n \to \infty} p_{\text{MD}}^{(n)} \left( R_{\boldsymbol{x}}^{(n)}; \theta_n' \right) \right] \leq \sup_{\{\theta_n'\} \in \tilde{\Theta}_{\{\theta_n\}}} \left[ \limsup_{n \to \infty} p_{\text{MD}}^{(n)} \left( \tilde{R}_{\boldsymbol{x}}^{(n)}; \theta_n' \right) \right],$$

simultaneously for all $\tilde{\Theta}_{\{\theta_n\}}$.

## V. AUMM Performance of Locally Asymptotic Normal Models

The results of section III can be extended, in the asymptotic sense defined in Section IV, to sequences of models which are LAN. An important example of LAN models is a memoryless sequence from a smooth parametric family. However, LAN models also include sequences with memory such as autoregressive processes, general ergodic Markov chains with smooth transition densities and Gaussian time series. We will survey some examples in the sequel, for more see, e.g., [17], [20], [21] and [22, Chapter 9]. Our goal is to establish sufficient conditions under which an AUMM rule can be found for such models. We start with giving the necessary background for LAN models followed by our results for AUMM without and with training data.

### A. Background

*1) LAN and the relation to the normal model:* A sequence of statistical models is LAN if, asymptotically, their likelihood-ratio processes are similar to those for the NLP [17]. The technical requirement is that the likelihood-ratio processes admit a certain quadratic expansion in a small neighborhood of $\theta_0$, as follows.

*Definition 13:* Suppose that $\theta_0$ is an inner point of $\Theta$ and that $\mathbf{Y} \sim \Pi_{\theta_0}^{(n)}$. Consider a sequence of models $\left( \Pi_\theta^{(n)} \Big| \theta \in \Theta \right)_{n \in \mathbb{N}}$. If there exists a sequence of invertible $[k \times k]$ matrices $r_n$ with

$r_n^{-1} \to 0$, a matrix $J_{\theta_0}$, and random vectors $u_n$ such that $u_n \rightsquigarrow \mathcal{N}_k(0, J_{\theta_0})$ under $\theta_0$, and for every converging sequence $h_n \to h$, inducing a sequence of parameters

$$\theta_n = \theta_0 + r_n^{-1} h_n, \tag{21}$$

the log likelihood-ratios satisfy

$$\log \frac{L_{\theta_n}(\mathbf{Y})}{L_{\theta_0}(\mathbf{Y})} = h^t u_n - \frac{1}{2} h^t J_{\theta_0} h + \delta_n, \tag{22}$$

where $\delta_n$ is a term that converges to zero in probability under $\Pi_{\theta_0}^{(n)}$, then we say that $\left(\Pi_\theta^{(n)} \middle| \theta \in \Theta\right)_{n \in \mathbb{N}}$ is LAN at $\theta_0$ with norming matrices $r_n$ and information matrix $J_{\theta_0}$.[11]

From here on we will assume that our sequence of models is indeed LAN at $\theta_0$ with norming matrices $r_n$ and with a positive definite information matrix $J_{\theta_0}$. We define the square-root of the information matrix $J_{\theta_0}^{1/2}$ as in (7).

Notice that since $\theta_0$ is an inner point, then for sufficiently large $n$, $\theta_n$ (21) is in $\Theta$. The measures $\Pi_{\theta_n}^{(n)}$ may be defined arbitrarily if $\theta_n \notin \Theta$. Typically the norming matrices are $r_n = \sqrt{n} I$ but not always. $r_n$ does not necessarily grow as $\sqrt{n}$, see examples 6 and 7 in Section V-C in the sequel. Also, an example of non-diagonal norming matrices can be found in [23]. In the univariate case ($k = 1$) the matrices $r_n$ are referred to as the norming rate. Under some regularity conditions the (probability-limit of the) Fisher-information matrix "per observation" (normalized by $r_n^{-2}$) qualifies as $J_{\theta_0}$.

To show the relation to the normal location problem we will parametrize the model locally around $\theta_0$ as follows. Consider any converging sequences $h_n \to h$ and $\theta_n \to \theta_0$ as in Definition 13. We define the local parameters $\mu_n$ as:

$$\begin{aligned}\mu_n &= J_{\theta_0}^{1/2} r_n (\theta_n - \theta_0) \\ &= J_{\theta_0}^{1/2} h_n.\end{aligned} \tag{23}$$

We denote the model $\Pi_\theta^{(n)}$ in terms of this local parameter as $\Pi_\mu^{(n)}$. The parameter $\mu = 0$ corresponds to $\theta_0$ whereas $\mu \neq 0$ corresponds to alternatives. According to (22) the "difficulty" of discriminating between $\theta_0$ and an alternative $\theta_n$ is determined by the quantity

$$d_n = \|\mu_n\|. \tag{24}$$

---

[11]We will omit the norming matrices or information matrix when not relevant.

We associate a sequence of parameters with $\mu_n$ converging to $\mu$, with the NLP with the same $\mu_1 = \mu$, and denote the limit of $d_n$ by $d$. Before we elaborate on this relation we start with the similarity in the log likelihood-ratios. It is clear from Definition 13 that if $\mu_n$ converges to $\mu$ then the asymptotic distribution of the log likelihood-ratio satisfies

$$\log \frac{L_{\theta_n}(\mathbf{Y})}{L_{\theta_0}(\mathbf{Y})} \rightsquigarrow \mathcal{N}\left(-\frac{1}{2}d^2, d^2\right), \text{ under } \theta_0. \tag{25}$$

On the other hand

$$\log \frac{L_{\theta_n}(\mathbf{Y})}{L_{\theta_0}(\mathbf{Y})} \rightsquigarrow \mathcal{N}\left(\frac{1}{2}d^2, d^2\right), \text{ under } \theta_n, \tag{26}$$

see [2, page 552]. Using these two properties we will show in Lemma 2 that the model $\Pi_\theta$ is asymptotically degraded at $\theta_0$ and characterize the asymptotic minimax sets. In other words, the asymptotic distribution of the log likelihood-ratio under $\theta_0$ and under $\theta_n$ is the same as that of the log likelihood-ratio in (8) with $\mu_1 = \mu$.

We will use the term "statistic" for a measurable map from $\mathbf{Y}$ which may depend on $\theta_0$ but does not depend on $\mu$. A randomized statistic $T$ based on the observation $Y$ is a measurable map $T = T(Y, C)$, such that $Y$ and $C$ are independent and $C \sim \text{unif}[0, 1]$. We note that in most cases we consider, randomization is not required. It turns out that the limit distributions of any sequence of statistics is necessarily the distribution of a (randomized) statistic in the model $\mathcal{N}_k(\mu, I)$, just as we saw above for the log likelihood-ratios. Thus every limiting distribution can be "represented" as a distribution of a statistic in the model $\mathcal{N}_k(\mu, I)$. For this reason this result is referred to as a "representation theorem" which is stated in proposition 2 in the sequel. In fact it can be extended to other (non-normal) models, but this is beyond the scope of our work. This result shows the profound relation to the model $\mathcal{N}_k(\mu, I)$, which applies to a wide range of problems. Indeed, such sequences of statistics may include discriminant rules or estimators.

The "limit model" $\mathcal{N}_k(\mu, I)$ provides an absolute standard for *what can be achieved* asymptotically by a sequence of tests or estimators, in the form of a "lower bound": No sequence of statistical procedures can be asymptotically better than the "best" procedure in the limit model [17]. This holds as long as by "best" we mean any optimality criterion which depends on the asymptotic distribution of statistics under various parameters, e.g. our AUMM criterion. For further reading and proofs see, [17, Chapters 7,9] and for this result in a hypothesis testing context see e.g. [17, Chapter 15] and [2, Chapter 13]. We formalize this important result in the following proposition.

*Proposition 2 (Representation Theorem, [17]):* Assume that $\left(\Pi_\theta^{(n)} \middle| \theta \in \Theta\right)_{n \in \mathbb{N}}$ is LAN at $\theta_0$ with non-singular information matrix $J_{\theta_0}$. Let $\{T_n\}$ be a sequence of statistics in the models $\left(\Pi_\mu^{(n)} \middle| \mu \in \mathbb{R}^k\right)_{n \in \mathbb{N}}$ such that $\{T_n\}$ converges in distribution under every $\mu$. Then there exists a randomized statistic $T$ in the (single-observation) NLP $\mathcal{N}_k(\mu, I)$ such that $T_n \rightsquigarrow T$ under any $\mu$.

*2) Estimation:* For a sequence of estimators (statistics) $\hat{\theta}_n$, let (cf. (23))

$$\hat{\mu}_n = J_{\theta_0}^{1/2} r_n(\hat{\theta}_n - \theta_0). \tag{27}$$

When it comes to estimating the underlying parameter, efficiency in the Fisher sense is defined as follows.

*Definition 14:* A sequence of estimators $\hat{\theta}_n$ is called efficient at $\theta_0$ if

$$\hat{\mu}_n \rightsquigarrow \mathcal{N}_k(0, I), \text{ under } \theta_0.$$

In the case of i.i.d. observations, this definition reduces to the classical Fisher efficiency (see [24, Chapter 11] and [25]). For $\theta \in \tilde{\Theta}$ we do not require efficiency, but only that the estimators are uniformly $r_n$-consistent. That is, with high probability the estimators are within a range of order $r_n^{-1}$ of the true parameter, whatever it may be.

*Definition 15:* A sequence of estimators $\hat{\theta}_n$ is uniformly $r_n$-consistent if for every $\varepsilon > 0$ there exist $M$ and $n_0$ such that

$$\sup_{\mu \in \mathbb{R}^k} P\left(\|\hat{\mu}_n - \mu\| > M; \mu\right) < \varepsilon,$$

for all $n > n_0$.

*B. Main Result*

We start by proving that the model $\Pi_\theta$ is asymptotically degraded at $\theta_0$, thus AUMM is relevant to it; then we prove the existence of an AUMM rule.

*Lemma 2:* Let $\left(\Pi_\theta^{(n)} \middle| \theta \in \Theta\right)_{n \in \mathbb{N}}$ be LAN at $\theta_0$ with norming matrices $r_n$ and a non-singular information matrix $J_{\theta_0}$. Then the model $\Pi_\theta$ is asymptotically degraded at $\theta_0$. Furthermore, the family of asymptotic minimax sets can be written as

$$\tilde{\Theta}[d] = \left\{ \{\theta_n\} \subseteq \tilde{\Theta} \middle| \liminf_{n \to \infty} d_n \geq d \right\}, \tag{28}$$

for all $d > 0$, where $d_n$ is defined in (24).

*Proof:* Fix some $p_{\text{FA}}, \bar{p}_{\text{MD}} \in (0,1)$. Let $\{\theta_n\}$ be a sequence of alternatives and suppose that $\limsup p^*_{\text{MD}}(p_{\text{FA}}; \theta_n) = \bar{p}_{\text{MD}}$. Let $\theta_{n_j}$ be a subsequence such that

$$\lim_{j \to \infty} p^*_{\text{MD}}(p_{\text{FA}}; \theta_{n_j}) = \bar{p}_{\text{MD}} \tag{29}$$

It is clear from Definition 13 that unless $d_{n_j}$ is bounded away from $\infty$, both error probabilities for the LRT tend to 0. Moreover, for (29) to hold we must have $d_{n_j} \longrightarrow d > 0$. Using (25) and (26) for $\theta_{n_j}$ we can obtain the tradeoff between $p_{\text{FA}}$ and $\bar{p}_{\text{MD}}$. Specifically, using the same arguments which led to (9) yields

$$Q^{-1}(p_{\text{FA}}) + Q^{-1}(\bar{p}_{\text{MD}}) = d.$$

It follows that $\bar{p}_{\text{MD}}$ is decreasing in $d$ for all $p_{\text{FA}}$. Let $\{\theta'_n\}$ be a sequence of alternatives with the corresponding $\tilde{d}_n$ defined as in (24). The sequence $\{\theta'_n\}$ is harder than $\{\theta_n\}$ if, and only if, $\liminf \tilde{d}_n = \tilde{d} < d$. Thus, as in (10), the family of asymptotic minimax sets can be written as (28) for all $d > 0$. Since for every two such sequences of alternatives either $\tilde{d} \leq d$ or $d \leq \tilde{d}$ it follows that the model is asymptotically degraded at $\theta_0$. ∎

Our main result is given in the next theorem. It shows how certain efficient estimators can be used to construct AUMM rules without and with training data. In both cases the AUMM performance can be described using the UMM performance in (18). The discriminant rule that we find asymptotically optimal is the following. For an LAN model, denote by $\hat{\theta}_{y,n}$ and $\hat{\theta}_{x,n_X}$ estimators of $\theta$ based on $\mathbf{Y}$ and $\mathbf{X}$ respectively. Let $\hat{\mu}_{y,n}$ be according to (27) and let

$$\hat{\mu}_{x,n} = J_{\theta_0}^{1/2} r_n (\hat{\theta}_{x,n_X} - \theta_0).$$

We will assume that

$$\lim_{n \to \infty} r_{n_X} r_n^{-1} = \sqrt{\rho} I, \tag{30}$$

for some $\rho \geq 0$. We denote as $\underline{R}_{\boldsymbol{x}}$ the discriminant rule in (17) with $x$ and $y$ replaced by $\hat{\mu}_{x,n}$ and $\hat{\mu}_{y,n}$ respectively and with $\rho$ in (30). That is,

$$\|\rho \hat{\mu}_{x,n} + \hat{\mu}_{y,n}\|^2 < Q^{-1}_{(k),\eta_0}(p_{\text{FA}}), \tag{31}$$

where $\eta_0 = \|\rho \hat{\mu}_{x,n}\|^2$.

*Theorem 2:* Let $\left(\Pi_\theta^{(n)} \big| \theta \in \Theta\right)_{n \in \mathbb{N}}$ be LAN at $\theta_0$ with norming matrices $r_n$ and a non-singular information matrix $J_{\theta_0}$. Suppose there exists a sequence of uniformly $r_n$-consistent estimators $\hat{\theta}_n$ which is efficient at $\theta_0$ and that (30) holds for some $\rho \geq 0$. Then the AUMM performance with guaranteed $p_{\text{FA}}$ is given by $\underline{p}_{\text{MD}}^{\text{umm}}(p_{\text{FA}}, d, \rho, k)$ in (18) for all minimax sets (28) with $d > 0$, and furthermore it is achieved by $\underline{R}_{\boldsymbol{x}}$ in (31).

The proof of Theorem 2 is immediate from the following two lemmas which yield the "direct" and "converse" parts.[12] They are proven in Appendix B. For these lemmas we assume that $\rho > 0$. The case $\rho = 0$ follows follows from the following consideration: On one hand the AUMM error probability must be a non-increasing function of $\rho$, thus

$$\underline{p}_{\text{MD}}^{\text{umm}}(p_{\text{FA}}, d, 0, k) \geq \lim_{\rho \downarrow 0} \underline{p}_{\text{MD}}^{\text{umm}}(p_{\text{FA}}, d, \rho, k),$$

and on the other hand for all $\rho$, $\underline{p}_{\text{MD}}^{\text{umm}}(p_{\text{FA}}, d, \rho, k)$ is at most the UMM $p_{\text{MD}}$ without training data (obtained by the GLRT), given by (12). As these two bounds coincide, the expression also extends to $\rho = 0$. Furthermore notice that for $\rho = 0$, $\underline{R}_{\boldsymbol{x}}$ ignores the training data since (31) does not depend on $\mathbf{X}$.

*Lemma 3 (Achievability):* Assume the conditions of Theorem 2 hold and that $\rho > 0$. Let $\{\theta_n\}$ be any sequence of alternatives with the corresponding $d_n$ as in (24). Then $\underline{R}_{\boldsymbol{x}}$ has a guaranteed significance level $p_{\text{FA}}$ and its asymptotic worst case missed-detection probability against $\{\theta_n\}$ is given by

$$\limsup_{n \to \infty} p_{\text{MD}}^{(n)}(\underline{R}_{\boldsymbol{x}}; \theta_n) = \underline{p}_{\text{MD}}^{\text{umm}}(p_{\text{FA}}, d, \rho, k), \tag{32}$$

where $d = \liminf d_n$.

The proof of this lemma uses (25) and Le Cam's third lemma [2, Corollary 12.3.2] to show that $\hat{\mu}_{y,n}$ and $\hat{\mu}_{x,n}$ are asymptotically distributed as $Y$ and $X$ in (17). Therefore the test statistic in (31) has the same asymptotic distribution as (17) with $\Delta = d$.

*Lemma 4 (Lower bound):* Assume the conditions of Theorem 2 hold and that $\rho > 0$. Let $R_{\boldsymbol{x}}^{(n)}$ be a sequence of rules for the problems in (20) with guaranteed significance level $p_{\text{FA}}$. Then the asymptotic worst case missed-detection probability satisfies

$$\sup_{\{\theta_n\} \in \tilde{\Theta}[d]} \left[ \limsup_{n \to \infty} p_{\text{MD}}^{(n)}\left(R_{\boldsymbol{x}}^{(n)}; \theta_n\right) \right] \geq \underline{p}_{\text{MD}}^{\text{umm}}(p_{\text{FA}}, d, \rho, k).$$

---

[12] We note that Lemma 4 generalizes [2, Theorem 13.5.4] which applies for an i.i.d. sequence from a smooth family of distributions without training data.

The proof of this lemma uses a version of Proposition 2 specialized for hypothesis testing, for the product model

$$\left\{ \Pi_{\tilde{\theta}}^{(n_X)} \times \Pi_\theta^{(n)} \,\middle|\, (\tilde{\theta}, \theta) \in \Theta \times \Theta \right\}$$

which is shown to be LAN at $(\theta_0, \theta_0)$ with norming matrices $\Psi_n$ and an information matrix $\Lambda$ given by the following $[2k \times 2k]$ block-matrices

$$\Psi_n = \begin{bmatrix} r_{n_X} & 0 \\ 0 & r_n \end{bmatrix}, \Lambda = \begin{bmatrix} \rho J_{\theta_0} & 0 \\ 0 & J_{\theta_0} \end{bmatrix}.$$

*C. Examples*

In this section we give examples for LAN models, thus we demonstrate the applicability of our main result. For all of these examples Theorem 2 applies, that is, the UMM performance is given by (18), up to a vanishing correction term $\delta_n$ which depends upon the specific example.

*1) i.i.d. models:* An important example of LAN models is repeated sampling from a smooth (differential in quadratic mean) distribution $P_\theta$. Specifically, let $\Pi_\theta^{(n)}$ be the model corresponding to $n$ i.i.d. observations from $P_\theta$. Then the sequence of models $\left( \Pi_\theta^{(n)} \,\middle|\, \theta \in \Theta \right)_{n \in \mathbb{N}}$ is LAN at every (inner point) $\theta$ with norming matrices $r_n = \sqrt{n} I$ and with information matrix $J_\theta$ which equals the Fisher-information matrix of a single observation from $P_\theta$. Example 1 in section IV is a special case with $P_\theta = \mathcal{N}_k(\theta, I)$. In the following two examples we show other common distributions.

*Example 2:* Suppose that $\Pi_\theta^{(n)}$ corresponds to $n$ i.i.d. observations from an exponential family with the following density

$$p_\theta(y) = g(\theta) h(y) e^{\eta(\theta)^t T(y)}.$$

Under some regularity conditions on the maps $\theta \mapsto \eta(\theta)$, the sequence of models $\left( \Pi_\theta^{(n)} \,\middle|\, \theta \in \Theta \right)_{n \in \mathbb{N}}$ is LAN with norming rate $\sqrt{n}$ and Fisher-information given by

$$J_\theta = \nabla \eta(\theta) \mathrm{Cov}_\theta(T(Y)) (\nabla \eta(\theta))^t,$$

see [17, page 96]. As estimators one can usually take the MLEs $\hat{\theta}_{y,n}$ and $\hat{\theta}_{x,n}$. The expressions in the LHS of (31) are given by

$$\begin{aligned} \hat{\mu}_{y,n} &= J_{\theta_0}^{1/2} \sqrt{n} (\hat{\theta}_{y,n} - \theta_0) \\ \hat{\mu}_{x,n} &= J_{\theta_0}^{1/2} \sqrt{n} (\hat{\theta}_{x,n_X} - \theta_0) \\ \rho &= \lim \frac{n_X}{n}. \end{aligned} \quad (33)$$

*Example 3:* Suppose that $\Pi_\theta^{(n)}$ corresponds to $n$ i.i.d. observations from a distribution on a finite alphabet $|\mathcal{Y}| = m$ defined by

$$P_\theta = \left(p_\theta^{(1)}, \ldots, p_\theta^{(m)}\right).$$

$P_\theta$ is parametrized by $\theta$ which holds the first $k = (m-1)$ probabilities. If $P_{\theta_0}$ has strictly positive probabilities then the sequence of models $\left(\Pi_\theta^{(n)} \middle| \theta \in \Theta\right)_{n \in \mathbb{N}}$ is LAN at $\theta_0$ with $r_n = \sqrt{n}I$. The Fisher information is the $[k \times k]$ matrix $J_{\theta_0} = (\tau_{i,j})$, where

$$\tau_{i,j} = \begin{cases} \frac{1}{p_{\theta_0}^{(i)}} + \frac{1}{p_{\theta_0}^{(m)}}, & i = j \\ \frac{1}{p_{\theta_0}^{(m)}}, & i \neq j. \end{cases}$$

As estimators one can take the MLEs which are given by the empirical probability functions (the "types") $P_Y$ and $P_X$. The expressions in the LHS of (31) are given by (33) with

$$\begin{aligned} \hat{\theta}_{y,n} &= P_Y, \\ \hat{\theta}_{x,n} &= P_X. \end{aligned} \tag{34}$$

Interestingly, for $\rho = 0$ the test statistic reduces to the well known Pearson $\mathcal{X}^2$-statistic. That is,

$$\|\hat{\mu}_{y,n}\|^2 = n\mathcal{X}^2\left(P_Y \| P_{\theta_0}\right), \tag{35}$$

where

$$\mathcal{X}^2\left(P_Y \| P_{\theta_0}\right) = \sum_{i=1}^m \frac{\left(p_Y^{(i)} - p_{\theta_0}^{(i)}\right)^2}{p_{\theta_0}^{(i)}}$$

is the $\mathcal{X}^2$-divergence. Furthermore, the quantity $d_n^2$ in (24) turns out to be

$$d_n^2 = n\Delta_n^2,$$

$$\Delta_n^2 = \mathcal{X}^2\left(P_{\theta_n} \| P_{\theta_0}\right).$$

As all the expressions we use are "local", the Chi-square divergence approximates up to a constant factor any smooth $f$-divergence [26, Sec. 4], and it can be thus replaced by e.g., the KL divergence. Since the estimators in (34) are sum of i.i.d. random variables and the acceptance region is a convex set for $\hat{\mu}_{y,n}$, we can identify the correction term in Theorem 2 for this case. Specifically, using lemma 14.4.1 in [2], which can be viewed as a multivariate version of the the Berry-Esseen Theorem, it turns out the correction term is non-uniform in $k$ and is given by $O\left(k^{1/4}/\sqrt{n}\right)$. We revisit this example in section VI.

*2) Wide-sense stationary models:* As mentioned, LAN models need not be memoryless. A basic example is a $k$-th order autoregressive (AR) model.

*Example 4:* Consider the following stable AR$(k)$ process parametrized by $\theta = \left(\theta^{(1)}, \ldots, \theta^{(k)}\right)$

$$Y_t = \theta^{(1)} Y_{t-1} + \theta^{(2)} Y_{t-2} + \ldots \theta^{(k)} Y_{t-k} + \varepsilon_t,$$

where the initial conditions are chosen to satisfy wide-sense stationarity, [13] and where $\{\varepsilon_t\}_{t=1}^n$ are i.i.d. $\mathcal{N}(0, \sigma_\varepsilon^2)$ for some known $\sigma_\varepsilon > 0$. Suppose that the roots of the "characteristic polynomial"

$$1 - \sum_{j=1}^{k} \theta^{(j)} z^j,$$

lie outside the unit circle. Let $\Pi_\theta^{(n)}$ correspond to the first $n$ observations from this model. Then the sequence of models $\left(\Pi_\theta^{(n)} \middle| \theta \in \Theta\right)_{n \in \mathbb{N}}$ is LAN at any $\theta$ with norming rate $r_n = \sqrt{n} I$ and with

$$J_\theta = \frac{1}{\sigma_\varepsilon^2} \Sigma_k(\theta),$$

where $\Sigma_k(\theta)$ is the auto-covariance matrix of $k$ successive samples of the process. As estimators one can take the MLEs $\hat{\theta}_{y,n}$ and $\hat{\theta}_{x,n}$. The expressions in the LHS of (31) are as in (33).

This result extends to other distributions for the noise which satisfy some regularity conditions on the Fisher-information. The observations in this model form a stationary Markov chain. The result also extends to general ergodic Markov chains with smooth transition densities (see [17, page 104] and [20]).

*3) Miscellaneous examples:* In this section we present several examples which show the following:

(i) Using the MLE may result in trivial performance.
(ii) There are LAN models with different norming matrices, even for i.i.d. models.
(iii) The condition of Theorem 2 on the norming matrices does not always hold.

We present the examples accordingly, starting with (i). Notice that the consistency of the MLE requires *global* "good behavior" whereas the LAN property is a *local* approximation. Therefore the LAN property does not imply the consistency of the MLE. In this example we show that using the GLRT yields trivial performance.

---

[13] Clearly, the asymptotic results hold also if the process is only asymptotically stationary.

*Example 5:* This example is based on [27]. Let $\Pi_\theta^{(n)}$ be the model corresponding to $n$ i.i.d. observations from the density $f_\theta$ on $[0,1]$ which is a mixture of Beta$(1,1)$ (uniform) and Beta$(\alpha, \beta)$

$$f_\theta(y) = \theta g(y; 1, 1) + (1-\theta) g(y; \alpha(\theta), \beta(\theta)),$$

where $g(y; \alpha, \beta)$ is the Beta$(\alpha, \beta)$ density given by

$$g(y; \alpha, \beta) = \frac{\Gamma(\alpha+\beta)}{\Gamma(\alpha)\Gamma(\beta)} y^{\alpha-1}(1-y)^{\beta-1} 1_{[0,1]}(y).$$

The parameter space is $\Theta = [\frac{1}{2}, 1]$ and $\alpha(\theta), \beta(\theta)$ are defined as follows

$$\alpha(\theta) = \theta \delta(\theta)$$
$$\beta(\theta) = (1-\theta)\delta(\theta)$$
$$\delta(\theta) = (1-\theta)^{-1} \exp\left((1-\theta)^{-2}\right).$$

For $\theta = 1$, $f_\theta$ is defined to be Beta$(1,1)$. In this example, The MLE for $\theta$ exists and converges to 1, *regardless of the true value of the parameter*, making the performance of the GLRT trivial. Nevertheless, Cramer's conditions ([28, page 500]) are satisfied, therefore there exists a sequence of asymptotically efficient estimators (which are still roots of the likelihood equation), $\hat{\theta}_y$, that are asymptotically normally distributed

$$\sqrt{n}(\hat{\theta}_y - \theta) \rightsquigarrow \mathcal{N}(0, J_\theta^{-1}),$$

when the true value $\theta$ is an interior point of $\Theta$. Take for example $\theta_0 = 2/3$ and consider the problem (2). The Fisher information is positive and finite in a neighborhood of $\theta_0$ and it is continuous in the parameter (since the beta distribution has a continuous Fisher-information matrix and $\delta(\theta)$ is twice continuously differentiable). The map $\theta \mapsto \sqrt{f_\theta(y)}$ is continuously differentiable for every $y$. Therefore, according to [17, lemma 7.6], the sequence of models $\left(\Pi_\theta^{(n)} \big| \theta \in \Theta\right)_{n \in \mathbb{N}}$ is LAN at $\theta_0$ with norming rate $\sqrt{n}$. It follows that the expressions in the LHS of (31) are as in (33) with $\rho = 0$.

In all of the examples so far the norming matrices were $r_n = \sqrt{n} I$. The following example shows that this is not always the case even for i.i.d. models.

*Example 6:* Suppose that $\Pi_\theta^{(n)}$ corresponds to $n$ i.i.d. observations from the density $f(y-\theta)$, where $f(y) = (1-|y|)^+$, is the triangular density. The sequence of models $\left(\Pi_\theta^{(n)} \big| \theta \in \Theta\right)_{n \in \mathbb{N}}$ is LAN at $\theta = 0$ with norming rate $r_n = \sqrt{n \log(n)}$. The existence of singularities in the density makes the estimation of the parameter easier, and hence a faster rescaling rate is necessary, see [17, pages 105, 212]. The parameter can be estimated using the MLE and the expressions in the LHS of (31) are as in (33).

We conclude this section with and example showing that the condition in Theorem 2 on the norming matrices does not always hold.

*Example 7:* Consider the following (non-stationary) linear-trend process $Y_j = \alpha + \delta \cdot j + \varepsilon_j$ parametrized by $\theta = (\alpha, \delta)$, where $\varepsilon_j \sim \mathcal{N}(0, \sigma^2)$, i.i.d., for some known $\sigma > 0$. Let $\Pi_\theta^{(n)}$ be the model corresponding to the first $n$ observations. Then $\left(\Pi_\theta^{(n)} \middle| \theta \in \Theta\right)_{n \in \mathbb{N}}$ is LAN at every $\theta$ with norming matrices and information matrix given by

$$r_n = \begin{bmatrix} \sqrt{n} & 0 \\ 0 & n^{\frac{3}{2}} \end{bmatrix}, J_\theta = \frac{1}{\sigma^2} \begin{bmatrix} 1 & \frac{1}{2} \\ \frac{1}{2} & \frac{1}{3} \end{bmatrix}.$$

For any $n_X$ and $n$ the norming matrices satisfy

$$r_{n_X} r_n^{-1} = \begin{bmatrix} \sqrt{\frac{n_X}{n}} & 0 \\ 0 & \left(\frac{n_X}{n}\right)^{3/2} \end{bmatrix},$$

which converges to a finite limit only when $(n_X/n) \to \tau$, for some scalar $\tau \geq 0$. Thus condition (30) does not hold for $\tau > 0$.

## VI. THE UMM BLOCKLENGTH-DIMENSION TRADEOFF

In this section we explore the asymptotic effect of the different parameters on the UMM performance. We consider sequences of problems with parameter values such that the asymptotic UMM performance does not degenerate, i.e., becomes trivial or perfect. We start with the simple case of a sequence of normal location problems with a single observation followed by normal location problems with blocklegths $n, n_X$. Finally we attempt to extend our results to finite alphabet sequences. As we shall see such an extension is not straightforward.

### A. Dimension and training asymptotics for the NLP

We investigate the UMM tradeoff expression (18). Consider the following sequence of parameters $\{(\Delta_m, \rho_m, k_m)\}_{m=1}^\infty$. Denote the UMM missed-detection probabilities by

$$p_{\text{MD}}^{(m)} = \underline{p}_{\text{MD}}^{\text{umm}}(p_{\text{FA}}, \Delta_m, \rho_m, k_m), \tag{36}$$

where $p_{\text{FA}} \in (0,1)$ is fixed. The following theorem gives an expression for the asymptotic behavior of $p_{\text{MD}}^{(m)}$ in terms of the UMM *hardness parameter*

$$\mathcal{E}_m = \frac{\Delta_m^2(1+2\rho_m)}{\sqrt{2k_m(1+2\rho_m)+4(1+\rho_m)^2\Delta_m^2}}. \tag{37}$$

*Theorem 3:* Suppose that $\{\mathcal{E}_m\}_{m=1}^\infty$ in (37) is bounded away from 0 and $\infty$. Then there exists a sequence

$$\psi_m = O\left(\frac{1}{\sqrt{\max\{k_m,\rho_m\}}}\right),$$

such that for any $p_{\text{FA}} \in (0,1)$, $p_{\text{MD}}^{(m)}$ in (36) there is a point $(p'_{\text{FA}}, p'_{\text{MD}})$ on the curve

$$Q^{-1}(p'_{\text{FA}}) + Q^{-1}(p'_{\text{MD}}) = \mathcal{E}_m \tag{38}$$

satisfying

$$\max\left\{|p_{\text{FA}} - p'_{\text{FA}}|, |p_{\text{MD}}^{(m)} - p'_{\text{MD}}|\right\} \leq \psi_m.$$

The proof is given in appendix C. It relies mainly on the normal approximation to the non-central $\mathcal{X}^2$ distribution. Notice that for any $p_{\text{FA}} \in (0,1)$, if $\mathcal{E}_m$ tends to zero or infinity then $p_{\text{MD}}^{(m)}$ becomes trivial or $p_{\text{MD}}^{(m)} \to 0$ respectively. The interpretation of the theorem is as follows. Consider the regime where the correction term is small (high $k$ or high $\rho$), then:

(i). The UMM hardness parameter $\mathcal{E}_m$ encapsulates the LRT hardness $\Delta$, the dimension $k$ and the amount of training $\rho$. If these factors change in a way that keeps $\mathcal{E}_m$ fixed, the asymptotic performance remains the same.

(ii). The sum of inverse Q-functions measures the UMM performance, similar to the requirement in (9). In this sense it is a good *single-letter* performance measure reflecting the entire $\left(p_{\text{FA}}, p_{\text{MD}}^{(m)}\right)$ curve.

In some cases, we can simplify the expression for $\mathcal{E}_m$. Specifically, we take the cases where one of the terms in the denominator is dominant. This is summarized in the following two corollaries, where we also restrict our attention to conditions under which $\mathcal{E}_m$ remains bounded away from zero and infinity, as in the theorem. The proofs are immediate, thus they are omitted.

*Corollary 1:* Suppose that $k_m = o(\rho_m)$ and that $\Delta_m = \Theta(1)$. Then Theorem 3 holds with $\psi_m = O\left(1/\sqrt{\rho_m}\right)$, and

$$\mathcal{E}_m = \Delta_m \left[1 - \frac{1}{2(1+\rho_m)}\right] + O\left(\frac{k_m}{\rho_m}\right).$$

*Corollary 2:* Suppose that $k_m = \Theta(\Delta_m^4(1+\rho_m))$ and that $\Delta_m = \Omega(1)$. Then Theorem 3 holds with

$$\mathcal{E}_m = \frac{\Delta_m^2 \sqrt{1+2\rho_m}}{\sqrt{2k_m}} + O\left(\sqrt{\frac{1+\rho_m}{k_m}}\right)$$

We interpret these expressions as follows. The tradeoff curve in (38) with $\mathcal{E}_m \to \Delta_m$ corresponds, asymptotically, to the performance of the LRT in (9). Of course, for this we will need a "sufficient" quality of training data compared to the dimension, as specified in Corollary 1. However, in practice we may not have such high quality training data. Then, as Corollary 2 shows, in order to compensate for the growing dimension we should have $\Delta_m$ growing fast enough.

## B. Introducing blocklength: Multi-measurement NLP

The UMM tradeoff (18) analyzed in Section VI-A gives the optimal performance also for the multiple-measurement version of the NLP, as discussed in Example 1. This allows us to introduce blocklengths $n_X$ and $n$ into the asymptotic picture. Some questions naturally arise: How should the amount of data grow, in order to compensate for the deteriorating performance due to the growing dimensionality $k$? If we are given a total fixed amount of data, how should we allocate it between training and test?

In addressing these questions, we index the problem using the blocklength $n$, rather than the index $m$ used above. Namely, consider a sequence of problems such as in example 1 indexed by $n$. Unlike the sequences of problems in section IV, we also allow the dimension $k$ to vary with $n$. To be more specific, for every $n$ we face the problem (20) corresponding to the sequence of parameters $n_X = n_X(n)$, $k_n$ and $\theta_0^{(k_n)} = 0$. Let $\theta_n^{(k_n)}$ be a sequence of alternatives. For every $n$, the UMM performance is given by

$$p_{\text{MD}}^{(n)} = \underline{p}_{\text{MD}}^{\text{umm}}(p_{\text{FA}}, d_n, \rho_n, k_n), \tag{39}$$

where

$$\begin{aligned}
d_n &= \sqrt{n}\Delta_n, \\
\Delta_n &= \|\theta_n^{(k_n)}\|, \\
\rho_n &= \frac{n_X}{n}.
\end{aligned} \tag{40}$$

Clearly corollaries 1 and 2 hold with (40) when replacing $\Delta_m$ with $d_n$.

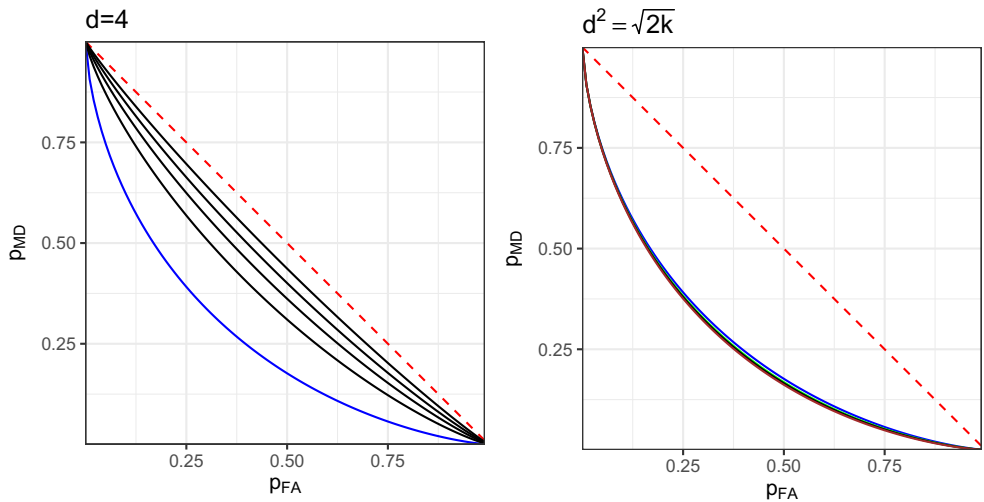

Fig. 4. The UMM performance in high dimension without training data, $p_{\mathrm{MD}} = \underline{p}_{\mathrm{MD}}^{\mathrm{umm}}(p_{\mathrm{FA}}, d, 0, k)$. The blue solid curve is the required tradeoff curve corresponding to $\mathcal{E} = 1$ and the red dashed line is the trivial performance. We take increasing dimension: $k = 128, 500, 1000, 2000, 5000$. In the left plot $d = \sqrt{n}\Delta$ is fixed at $d = 4$ therefore the performance becomes trivial as $k$ grows. In the right plot $d^2 = \sqrt{2k}$, maintaining the same tradeoff curve corresponding to $\mathcal{E} = 1$.

*1) Blocklength-dimension tradeoff:* Keeping $\rho_n$ and $\Delta_n$ fixed, we can consider the tradeoff between $n$ and $k_n$. Using Corollary 2, we see that in order to keep the performance fixed, $k_n$ grows as $n^2$. In other words, for a sequence of problems with growing dimension $k_n$ satisfying the conditions of Corollary 2, the blocklength $n$ must grow as $\sqrt{k_n}$ to avoid trivial performance. Figure 4 demonstrates the effect for $\rho_n = 0$ and $\mathcal{E} = 1$.

*2) Optimal allocation of data between test and training:* Keeping $\Delta_n$ and $k_n$ fixed, we can ask the following, Suppose that we have a total sample budget allows a fixed number of samples $n_T$ that can be divided between training blocklength $n_X$ and test blocklength $n$ as we wish, that is

$$n_X + n = n_T.$$

What policy optimizes performance? In the limit of interest, we can answer this asymptotically by maximizing $\mathcal{E}_m$ in (37) over $\rho$, when

$$a_n = n_T \Delta^2 \tag{41}$$
$$= (1 + \rho_n) n \Delta^2$$

is held fixed at some $a \in \mathbb{R}$. Specifically, for the sequence of problems with (40), by (39), Theorem 3 holds. The corresponding expression in (37) with (41) becomes

$$\begin{aligned}\mathcal{E}_n(\rho_n) &= \frac{n\Delta^2(1+2\rho_n)}{\sqrt{2k(1+2\rho_n)+4(1+\rho_n)^2 n\Delta^2}} \\ &= \frac{a(1+2\rho_n)}{(1+\rho_n)\sqrt{2k(1+2\rho_n)+4(1+\rho_n)a}}.\end{aligned} \quad (42)$$

By differentiating the RHS of (42) w.r.t. $\rho_n$, one finds that when $k_n$ is large the (asymptotically) optimal allocation will be achieved for $\rho_n = 0$. In other words, in high-dimensional NLP, all data should be test data.

The conclusion may seem counterintuitive to a reader that is familiar with machine-learning literature. However, we point out two main differences between the problem at hand and a typical supervised-learning classification scenario. First, out of the two classes we have, only one needs to be learned while the statistics of the other is a-priori known. Second, and perhaps more importantly, the growing test data is guaranteed to entirely belong to one of the classes, thus as the size grows the "synergy" between the test samples improves.

## C. Extension to LAN models

Naturally, we would like to extend the results of Section VI-A beyond the normal model. However, caution is needed. If we keep the dimension $k_n$ fixed, then by Theorem 2 the UMM performance is given by (39), with a correction term that vanishes as $n$ grows. However, one may be interested in a family of LAN models indexed by the dimension, such that for each blocklength $n$ we have a member of a model with the corresponding $k_n$. Unfortunately, in general the convergence of the correction to (39) is not uniform in $k_n$. Thus, if $k_n$ grows "too fast" with respect to $n$, we may not use the conclusions drawn from Theorems 2 and 3.

*Example 8:* Consider the sequence of discrete i.i.d. problems of Example 3, indexed by $n$. Fix $\rho = 0$, and for any dimension $k_n$ let $\theta_0^{(k_n)}$ correspond to the uniform distribution over $m = k_n + 1$ coordinates. It is shown in [29] that the test statistic (35) converges to the $\mathcal{X}^2$-distribution if and only if $k_n = o(n^2)$. When this does not hold, our AUMM results cannot be used. Recalling that for fixed performance in the multi-measurement NLP problem we required $k_n = \Theta(n^2)$, obviously the conclusion cannot be applied to this example (but any slower growth of $k_n$ would do).

We conclude, that for the dimension-asymptotic analysis to be applied, the dimension cannot grow too fast. For any specific LAN model, one can analyze convergence and find conditions for convergence, as well as characterize the correction term.

*D. Discussion: Average significance Level*

In Section II-D we defined UMM detectors with guaranteed significance level. Indeed, we can also define UMM under average significance level as follows. The average false-alarm probability is given by

$$p_{\text{FA}}(R_x; \theta_1) = E\left[p_{\text{FA}}(R_x|X)\right] = P(R_X(Y) = 1; \theta_0, \theta_1)$$

and we require that for all minimax sets $\tilde{\Theta}_\theta$ and for all $p_{\text{FA}}$, the detector will minimize

$$\max_{\theta_1 \in \tilde{\Theta}_\theta} p_{\text{MD}}(R_x; \theta_1)$$

among all detectors with

$$\max_{\theta_1 \in \tilde{\Theta}_\theta} p_{\text{FA}}(R_x; \theta_1) \leq p_{\text{FA}}.$$

We argue that such detectors do not exist for the NLP. However in an asymptotic sense, when the UMM performance with guaranteed significance level converges to the curve (38), it also gives the *asymptotic* UMM performance with average significance level.

In order to see why strict UMM optimality with average significance level is not feasible, we revisit the detectors of Section III-B. Consider the minimax set (10) with $d = 1$, and assume that $\mu_1$ is uniform over the $d$-sphere, as in Lemma 1. Recall that the measurement $x$ induces an a-posteriori $\mu_1$ such that $\mu_1/\Delta$ has Von Mises-Fisher distribution with concentration parameter $(\rho\Delta\|x\|)$. The larger the concentration parameter is, the greater is the clustering around the mean direction, with zero corresponding to a uniform distribution, see [19, page 432]. Thus, the pair $(x,\|\Delta\|)$ induces a conditional trade-off curve. Now consider the minimax detector for the corresponding minimax set. It may still assume that $\mu_1$ is on the sphere (this is the worst-case assumption). Knowing $x$ and $\Delta$, it may choose a different working point $p_{\text{FA}}$ for each $\|x\|$-induced conditional curve, in a manner that minimizes the average (this will be achieved by keeping the slope of the conditional curve fixed for all $\|x\|$), as opposed to the guaranteed-significance detector which has to choose the same $p_{\text{FA}}$ for all $\|x\|$. Thus, the optimal detectors for different minimax sets do not coincide, and a UMM detector does not exist.

However, the curve (38) does give the asymptotic UMM performance with average significance level, whenever the conditions of Theorem 3 hold. This is due to the fact that when either $\rho_m$ or $k_m$ grow the variability of $x$ becomes inconsequential, thus approximately the same significance level is guaranteed for any $x$. To be more specific, denote the discriminant rule (15) by $R_x^*$ and suppose that the (average) false-alarm probability is $p_{\text{FA}}(R_x^*; \mu_1) = p_{\text{FA}}$, for some $p_{\text{FA}} \in (0,1)$. Recall from the proof of Theorem 1 in appendix A that the acceptance region of $R_x^*$ is also a sphere for $y$ centered at $(-\rho_m x)$ with radius involving the function $c_{k_m}(\cdot)$ and depending on $x$ only through $\sqrt{\vartheta_0} = \|\rho_m x\|$. Note that $\sqrt{\vartheta_0/\rho_m}$ has a non-central $\mathcal{X}$-distribution with $k$ degrees of freedom and non-centrality parameter $\sqrt{\rho_m}\Delta_m$. A known property of the $\mathcal{X}$-distribution is that as either $k$ or the non-centrality parameter grow, it converges to a normal distribution with a standard deviation which is negligible compared to its mean. Due to the smoothness and monotonicity of $c_{k_m}(\cdot)$, it can be shown that with probability which tends to 1 (w.r.t. $x$), the conditional false-alarm probability $p_{\text{FA}}(R_x^* \mid X)$, will be "close" to $p_{\text{FA}}$.

From an operational point of view, either guaranteed or average significance level may be of interest, according to the circumstances. In a setting where a training sequence is followed by a single test sequence, average significance level is more adequate. However imagine that a training sequence is available one time when the detector is designed, and then the same detector is used many times to test many sequences, waiting for an outlier. In that case, we wouldn't want a single "bad" training sequence to hurt performance, and therefore a guarantee on the significance level is desirable, and UMM detectors exist. The asymptotic view provides further justification for our UMM analysis: In the limit of high dimension or high training-data quality, it holds also when one is interested in average significance.

ACKNOWLEDGEMENT

The authors thank Meir Feder, Or Ordentlich and Avraham Sidi for helpful discussions.

# APPENDIX A
## PROOFS FOR SECTION III

### A. Proof of Lemma 1

The likelihoods are given by

$$L_{\mathcal{H}_0}(x,y) = L(y;0,I) \int_\Delta L\left(x;m,\frac{1}{\rho}I\right) dm$$

$$L_{\mathcal{H}_1}(x,y) = \int_\Delta L(y;m,I) \cdot L\left(x;m,\frac{1}{\rho}I\right) dm,$$

where $L(\cdot;\mu,\Sigma)$ is defined as in (8) and where the notation $\int_\Delta f(m)dm$ means the Lebesgue integral w.r.t. the uniform measure on the $\Delta$-sphere in $\mathbb{R}^k$. The likelihood-ratio satisfies

$$\frac{L_{\mathcal{H}_1}(x,y)}{L_{\mathcal{H}_0}(x,y)} = \frac{\hat{L}_{H_1}(y)}{L(y;0,I)},$$

where

$$\hat{L}_{H_1}(y) = \int_\Delta u(m) L(y;m,I) dm$$

$$u(m) = \frac{L\left(x;m,\frac{1}{\rho}I\right)}{\int_\Delta L\left(x;\tilde{m},\frac{1}{\rho}I\right) d\tilde{m}}.$$

For the LRT, the discriminant rule is $\mathbf{R}(y) = 0$ when

$$\log(\hat{L}_{H_1}(y)) - \log(L(y;0,I)) < T. \tag{43}$$

Notice that we can write

$$L(y;m,I) = [\det(2\pi I)]^{-\frac{1}{2}} \exp\left(-\frac{1}{2}\|y-m\|^2\right)$$

$$= [\det(2\pi I)]^{-\frac{1}{2}} \exp\left(-\frac{1}{2}\left[\|y\|^2 + \Delta^2\right]\right) \exp\left(y^t m\right),$$

and

$$L\left(x;m,\frac{1}{\rho}I\right) = \left[\det\left(\frac{2\pi}{\rho}I\right)\right]^{-\frac{1}{2}} \exp\left(-\frac{\rho}{2}\|x-m\|^2\right)$$

$$= \left[\det\left(\frac{2\pi}{\rho}I\right)\right]^{-\frac{1}{2}} \exp\left(-\frac{\rho}{2}\left[\|x\|^2 + \Delta^2\right]\right) \exp\left(\rho x^t m\right).$$

Therefore we have

$$L(y;m,I) \cdot L\left(x;m,\frac{1}{\rho}I\right) = [\det(2\pi I)]^{-\frac{1}{2}} \cdot \left[\det\left(\frac{2\pi}{\rho}I\right)\right]^{-\frac{1}{2}}$$

$$\cdot \exp\left(-\frac{1}{2}\left[\|y\|^2 + \rho\|x\|^2 + \Delta^2(1+\rho)\right]\right) \exp\left(m^t[\rho x + y]\right). \tag{44}$$

To evaluate $\hat{L}_{H_1}(y)$ we integrate the last term w.r.t. $\Delta$. We have

$$\int_\Delta \exp\left(m^t(\rho x + y)\right) dm = \frac{1}{c_k(\xi)} \int_1 c_k(\xi) \exp\left(\xi \tilde{m}^t w\right) d\tilde{m},$$

where $\tilde{m} = \frac{1}{\Delta}m$, $\xi = \Delta\|\rho x + y\|$ and $w = \frac{\Delta}{\xi}(\rho x + y)$ and with $c_k(\tau)$ as in (16). The last integrand is the p.d.f. of the von Mises-Fisher distribution with mean direction $w$ and concentration parameter $\xi$, see e.g. [19, Equation 15.3.6]. It follows that

$$\int_\Delta \exp\left(m^t[\rho x + y]\right) dm = \frac{1}{c_k(\xi)}. \tag{45}$$

Combining (45) with (44) and applying the logarithm function yields

$$\log\left(\int_\Delta L(y; m, I) \cdot L\left(x; m, \frac{1}{\rho}I\right) dm\right) = -\frac{1}{2}\log\left(\det(2\pi I)\right) - \frac{1}{2}\log\left(\det\left(\frac{2\pi}{\rho}I\right)\right)$$
$$- \frac{1}{2}\left[\|y\|^2 + \rho\|x\|^2 + \Delta^2(1+\rho)\right] - \log c_k(\xi).$$

Using similar arguments it can be shown that

$$\log\left(\int_\Delta L\left(x; m, \frac{1}{\rho}I\right) dm\right) = -\frac{1}{2}\log\left(\det\left(\frac{2\pi}{\rho}I\right)\right) - \frac{\rho}{2}\left[\|x\|^2 + \Delta^2\right] - \log\left(c_k\left(\Delta\|\rho x\|\right)\right).$$

This shows that

$$\log\left(\hat{L}_{H_1}(y)\right) = -\frac{1}{2}\log\left(\det(2\pi I)\right) - \frac{1}{2}\|y\|^2 + \log\left(\frac{c_k(\Delta\|\rho x\|)}{c_k(\Delta\|\rho x + y\|)}\right) - \frac{\Delta^2}{2}.$$

Finally, returning to (43), the LRT is equivalent to

$$\log\left(\frac{c_k(\Delta\|\rho x\|)}{c_k(\Delta\|\rho x + y\|)}\right) - \frac{\Delta^2}{2} < T,$$

which reduces to the required test.

## B. Proof of Theorem 1

The function $c_k(\cdot)$ is monotonously decreasing, see Lemma 5 below. Thus we can rewrite (15) as:

$$\|\rho x + y\|^2 < \Psi^2(\vartheta_0), \tag{46}$$

where

$$\Psi(\vartheta_0) = \frac{c_k^{-1}\left(T \cdot c_k(\Delta\sqrt{\vartheta_0})\right)}{\Delta}.$$

The conditional distribution of the LHS of (46) given the training sequence is given by

$$\|\rho X + Y\|^2 \mid X = x \sim \begin{cases} \mathcal{X}^2_{(k),\vartheta_0}, & \text{Under } \mathcal{H}_0 \\ \mathcal{X}^2_{(k),\vartheta_1}, & \text{Under } \mathcal{H}_1. \end{cases} \quad (47)$$

Therefore the rule (17) guarantees a significance level of $p_{\text{FA}}$. By (47) it is straightforward that the performance is given by (18). Although not immediately seen, the RHS of (18) depends upon the parameters only through $\Delta$ and $\rho$. To see this, notice that the distribution of the non-centrality parameters depends only on $\Delta$ and $\rho$. Specifically

$$\frac{\vartheta_i}{\rho} \sim \mathcal{X}^2_{(k),\lambda_i}, \ i = 0, 1 \quad (48)$$

where

$$\begin{aligned} \lambda_0 &= \rho \Delta^2, \\ \lambda_1 &= \frac{(1+\rho)^2}{\rho} \Delta^2. \end{aligned} \quad (49)$$

To prove that this is indeed the UMM performance, suppose to the contrary that (17) is not a UMM rule. Then there exist some $d > 0$, $p'_{\text{FA}} \in (0,1)$ and a different discriminant rule, $\tilde{R}_{\boldsymbol{x}}$ with guaranteed significance level $p'_{\text{FA}}$ such that

$$\sup_{\mu \in \tilde{\Theta}(d)} p_{\text{MD}}\left(\tilde{R}_{\boldsymbol{x}}; \mu\right) < \sup_{\mu \in \tilde{\Theta}(d)} p_{\text{MD}}\left(\underline{R}_{\boldsymbol{x}}; \mu\right),$$

where $\underline{R}_{\boldsymbol{x}}$ is the discriminant rule in (17) with $p_{\text{FA}}$ replaced by $p'_{\text{FA}}$. It can be shown that $\underline{p}^{\text{umm}}_{\text{MD}}(p_{\text{FA}}, \Delta, \rho, k)$ is a decreasing function of $\Delta$. Using (10), it follows that the supremum on the RHS is attained on the $d$-sphere. By (48) and (49), $\underline{p}^{\text{umm}}_{\text{MD}}(p_{\text{FA}}, \Delta, \rho, k)$ is constant on the $d$-sphere. It follows that

$$\sup_{\mu \in \tilde{\Theta}(d)} p_{\text{MD}}\left(\tilde{R}_{\boldsymbol{x}}; \mu\right) < \underline{p}^{\text{umm}}_{\text{MD}}(p'_{\text{FA}}, d, \rho, k). \quad (50)$$

In particular, if we limit the alternatives in the supermum in (50) to the $d$-sphere we get

$$p_{\text{MD}}\left(\tilde{R}_{\boldsymbol{x}}; \mu\right) < \underline{p}^{\text{umm}}_{\text{MD}}(p'_{\text{FA}}, d, \rho, k), \quad \forall \mu \text{ s.t. } \|\mu\| = d. \quad (51)$$

Therefore (51) also holds on average w.r.t. $\mu$ such that $\|\mu\| = d$. It follows that

$$p_{\text{MD}}\left(\tilde{R}_{\boldsymbol{x}}; d\right) < p_{\text{MD}}\left(\underline{R}_{\boldsymbol{x}}; d\right), \quad (52)$$

where $p_{\text{MD}}\left(\tilde{R}_{\boldsymbol{x}}; d\right)$ and $p_{\text{MD}}\left(\underline{R}_{\boldsymbol{x}}; d\right)$ are the missed-detection probabilities of $\tilde{R}_{\boldsymbol{x}}$ and $\underline{R}_{\boldsymbol{x}}$ for the problem (14) with $M$ uniformly distributed on the $d$-sphere, respectively. Since (52) holds on average w.r.t. $X$, it must be that for some value $\boldsymbol{x}$,

$$p_{\text{MD}}\left(\tilde{R}_{\boldsymbol{x}} \big| \boldsymbol{x}; d\right) < p_{\text{MD}}\left(\underline{R}_{\boldsymbol{x}} \mid \boldsymbol{x}; d\right). \tag{53}$$

Finally, since the distribution of $X$ is the same under both hypotheses, the likelihood-ratio for $Y \mid X$ and for $(X, Y)$ is the same. Hence the LRT for the problem (14) yields an optimal tradeoff between $p'_{\text{FA}}$ and $p_{\text{MD}}\left(\underline{R}_{\boldsymbol{x}} \mid \boldsymbol{x}; d\right)$. Since $\tilde{R}$ guarantees $p'_{\text{FA}}$, (53) contradicts the Neyman-Pearson lemma.

### C. Monotonicity of $c_k$

*Lemma 5:* The function $c_k(\cdot)$ is monotonously decreasing.

*Proof:* We will show that the derivative of $c_k(\cdot)$ is strictly negative. To this end we use the following known recurrence relation of the modified Bessel function

$$\tau \cdot \frac{d}{d\tau} \mathcal{I}_{(s)}(\tau) = \tau \mathcal{I}_{(s+1)}(\tau) + s \mathcal{I}_{(s)}(\tau),$$

for $s = (k/2 - 1)$, see [19, Equation 15.3.11]. Straightforward manipulation yields

$$\frac{d}{d\tau} c_k(\tau) = -\frac{\mathcal{I}_{(k/2)}(\tau)}{\mathcal{I}_{(k/2-1)}(\tau)} c_k(\tau).$$

Since the modified Bessel function is positive we get that the derivative of $c_k(\tau)$ is negative. ∎

## APPENDIX B
## PROOFS FOR SECTION V

### A. Proof of lemma 3

To prove Lemma 3, first we prove the following lemma.

*Lemma 6:* Assume the conditions of Theorem 2 hold with $\rho > 0$ and let $\underline{R}_{\boldsymbol{x}}$ be the discriminant rule in (31). Let $\{\theta_n\}$ be a sequence of alternatives with the corresponding $d_n$ as in (24). If $d_n \to \infty$ then

$$\lim_{n \to \infty} p_{\text{MD}}^{(n)}(\underline{R}_{\boldsymbol{x}}; \theta_n) = 0. \tag{54}$$

*Proof:* We can write

$$\hat{\mu}_{y,n} - \mu_n = J_{\theta_0}^{1/2} r_n(\hat{\theta}_{y,n} - \theta_n) + J_{\theta_0}^{1/2} r_n(\theta_n - \theta_0) - \mu_n \tag{55}$$
$$= J_{\theta_0}^{1/2} r_n(\hat{\theta}_{y,n} - \theta_n).$$

By assumption $\hat{\theta}_{y,n}$ is uniformly $r_n$-consistent, so the last term is bounded in probability. Therefore

$$\hat{\mu}_{y,n} - \mu_n = O_p(1), \text{ under } \theta_n.$$

Using similar arguments for $\hat{\mu}_{x,n}$ and the assumption that (30) holds, we have

$$\begin{aligned}\hat{\mu}_{x,n} - \mu_n &= J_{\theta_0}^{1/2} r_n(\hat{\theta}_{x,n_X} - \theta_0) - \mu_n \\ &= J_{\theta_0}^{1/2} r_n(\hat{\theta}_{x,n_X} - \theta_n) + J_{\theta_0}^{1/2} r_n(\theta_n - \theta_0) - \mu_n \\ &= J_{\theta_0}^{1/2} r_n(\hat{\theta}_{x,n_X} - \theta_n) \\ &= J_{\theta_0}^{1/2} \left(\frac{1}{\sqrt{\rho}} r_{n_X} + A_n\right)(\hat{\theta}_{x,n_X} - \theta_n),\end{aligned} \tag{56}$$

where $\{A_n\}$ is a sequence of matrices which converge to 0 element-wise. Since $\hat{\theta}_{x,n_X}$ is uniformly $r_{n_X}$-consistent the last term is bounded in probability, thus

$$\hat{\mu}_{x,n} - \mu_n = O_p(1), \text{ under } \theta_n,$$

and $\eta_0$ satisfies

$$\eta_0 = \rho^2 d_n^2 + O_p(1), \text{ under } \theta_n. \tag{57}$$

It follows that the test statistic in the LHS of (31) satisfies

$$\|\rho\hat{\mu}_{x,n_X} + \hat{\mu}_{y,n}\|^2 = (1+\rho)^2 d_n^2 + O_p(1), \text{ under } \theta_n, \tag{58}$$

where $d_n \to \infty$. Finally, to show that the asymptotic missed-detection probability of $\underline{R}_{\boldsymbol{x}}$ tends to 0 we use the following normal approximation to the non-central $\mathcal{X}^2$ distribution:

$$Q_{(k),\eta_0}^{-1}(p_{\text{FA}}) = \sqrt{2(k+2\eta_0)} Q^{-1}(p_{\text{FA}}) + k + \eta_0 + O\left(\frac{1}{\sqrt{\max\{k,\eta_0\}}}\right), \tag{59}$$

locally for any fixed $p_{\text{FA}} \in (0,1)$. The correction term on the RHS of (59) is $O\left(1/\sqrt{k}\right)$ due to Berry-Esseen but it is also $O\left(1/\sqrt{\eta_0}\right)$, see [30, page 466]. The significant tern in (59) is $\eta_0$. Using (57) and (58) we deduce that (54) holds.

Using Lemma 6, we prove Lemma 3 as follows.

*Proof:* First we show that $\underline{R}_{\boldsymbol{x}}$ has a guaranteed asymptotic significance level $p_{\text{FA}}$. Indeed, due to the efficiency of $\hat{\theta}_{y,n}$ at $\theta_0$ and the continuous mapping theorem, the asymptotic distribution of the test statistic given $\mathbf{X}$ under the null hypothesis is given by

$$\|\rho\hat{\mu}_{x,n} + \hat{\mu}_{y,n}\|^2 \mid \mathbf{X} \rightsquigarrow \mathcal{X}^2_{(k),\eta_0}, \text{ under } \theta_0.$$

Since the threshold is set by the RHS of (31), it follows that for every $\boldsymbol{x}$

$$\lim_{n\to\infty} p_{\text{FA}}^{(n)}(\underline{R}_{\boldsymbol{x}} \mid \boldsymbol{x}) = p_{\text{FA}}.$$

Denote

$$p = \limsup_{n\to\infty} p_{\text{MD}}^{(n)}(\underline{R}_{\boldsymbol{x}}; \theta_n).$$

If $d = \infty$ then $d_n \to \infty$ and according to Lemma 6, $p = 0$ which also equals $\underline{p}_{\text{MD}}^{\text{umm}}(p_{\text{FA}}, d, \rho, k)$. Thus we can ignore any subsequence of $\{d_n\}$ which tends to $\infty$, remaining only with bounded subsequences. Assume w.l.o.g. that $d_n$ itself is bounded. Specifically, assume that for all $n$: $\underline{d} \leq d_n \leq \overline{d} < \infty$.

We start with a simpler case where for all $n: \mu_n = \mu$, for some $0 \neq \mu \in \mathbb{R}^k$. Then $d = \|\mu\|$. The sequence $\hat{\theta}_{y,n}$ is efficient at $\theta_0$, therefore

$$\hat{\mu}_{y,n} \rightsquigarrow \mathcal{N}_k(0, I), \text{ under } \theta_0.$$

Using (25) and Le Cam's third lemma [2, Corollary 12.3.2], $\hat{\mu}_{y,n}$ is also asymptotically normal under $\theta_n$ with the same variance. Specifically, using (55) in the proof of Lemma 6 we have

$$\hat{\mu}_{y,n} \rightsquigarrow \mathcal{N}_k(\mu, I), \text{ under } \theta_n. \tag{60}$$

Furthermore, using similar arguments to those used for (56) together with Slutsky's lemma yields

$$\hat{\mu}_{x,n} \rightsquigarrow \mathcal{N}_k\left(\mu, \frac{1}{\rho}I\right), \text{ under } \theta_n. \tag{61}$$

Therefore the test statistic in (31) has the same asymptotic distribution as (17) with $\Delta = d$. It follows that

$$\begin{aligned} p &= \lim_{n\to\infty} p_{\text{MD}}^{(n)}(\underline{R}_{\boldsymbol{x}}; \theta_n) \\ &= \underline{p}_{\text{MD}}^{\text{umm}}(p_{\text{FA}}, d, \rho, k). \end{aligned} \tag{62}$$

Using the continuous mapping theorem for (60) and (61) it follows that (62) still holds if we replace $\mu$ with a converging sequence $\mu_n \to \mu$. Finally, suppose that $\{\mu_n\}$ is non-converging. Suppose $\{\theta_{n_j}\}$ is a subsequence of $\{\theta_n\}$ such that

$$\lim_{j \to \infty} p_{\text{MD}}^{(n_j)}(R_{\boldsymbol{x}}; \theta_{n_j}) = p.$$

Since $d_n \in [\underline{d}, \overline{d}]$ for all $n$, the subsequence $\{\mu_{n_j}\}$ lies in a compact set and it has a converging further subsequence. To simplify the notation assume w.l.o.g. that $\{\mu_{n_j}\}$ itself converges $\mu_{n_j} \to \mu$, where $\|\mu\| = d' \in [\underline{d}, \overline{d}]$. Thus according to (62) we have

$$p = \lim_{j \to \infty} p_{\text{MD}}^{(n_j)}(R_{\boldsymbol{x}}; \theta_{n_j})$$
$$= \underline{p}_{\text{MD}}^{\text{umm}}(p_{\text{FA}}, d', \rho, k).$$

Since $\underline{p}_{\text{MD}}^{\text{umm}}(p_{\text{FA}}, \cdot, \rho, k)$ is monotonically decreasing it must be that $d' = \liminf d_n$. Thus (32) holds for any $\rho > 0$. ∎

### B. Proof of Lemma 4

*Proof:* It suffices to find a sequence of alternatives $\{\theta_n\} \in \tilde{\Theta}[d]$ such that

$$\limsup_{n \to \infty} p_{\text{MD}}^{(n)}(R_{\boldsymbol{x}}^{(n)}; \theta_n) \geq \underline{p}_{\text{MD}}^{\text{umm}}(p_{\text{FA}}, d, \rho, k).$$

Fix some $0 \neq \mu \in \mathbb{R}^k$ and let $\theta_n$ be the alternative defined by the local parameter as in (23) with $\mu_n = \mu$ and $h_n = h$ for all $n$. That is

$$\begin{aligned} h &= J_{\theta_0}^{-1/2} \mu \\ &= r_n(\theta_n - \theta_0). \end{aligned} \quad (63)$$

With this local parametrization, denote

$$p = \limsup_{n \to \infty} p_{\text{MD}}^{(n)}\left(R_{\boldsymbol{x}}^{(n)}; h\right)$$

and suppose that $\{n_j\}$ is a sub-sequence such that

$$p_{\text{MD}}^{(n_j)}\left(R_{\boldsymbol{x}}^{(n_j)}; h\right) \longrightarrow p.$$

Consider the following product model

$$\left\{ \Pi_{\tilde{\theta}}^{(n_X)} \times \Pi_{\theta}^{(n)} \middle| (\tilde{\theta}, \theta) \in \Theta \times \Theta \right\} \quad (64)$$

with the local parametrizartion $(\tilde{h}, h)$ as in (63). According to Lemma 7 below the model in (64) is LAN at $(\theta_0, \theta_0)$ with norming matrices $\Psi_n$ and an information matrix $\Lambda$ given by the following $[2k \times 2k]$ block-matrices

$$\Psi_n = \begin{bmatrix} r_{n_X} & 0 \\ 0 & r_n \end{bmatrix}, \Lambda = \begin{bmatrix} \rho J_{\theta_0} & 0 \\ 0 & J_{\theta_0} \end{bmatrix}.$$

For the current proof we will use a version of Proposition 2 specialized for hypothesis testing given in [2, Theorem 13.4.1]. A discriminant rule $R$ is viewed as a statistic taking values in $[0, 1]$ with the interpretation that if $y$ is observed then a null hypothesis is rejected with probability $R(y)$. The null and alternative hypotheses will correspond to parameters of the form $(\tilde{h}, 0)$ and $(\tilde{h}, h), h \neq 0$, respectively. To be more specific, in the model (64) we denote the false-alarm and missed-detection probabilities of $R_{\boldsymbol{x}}^{(n)}$ with respect to the HT problem: $H_0 : (\tilde{h}, 0)$ vs. $H_1 : (\tilde{h}, h)$, as $p_{\text{FA}}^{(n)}\left(R_{\boldsymbol{x}}^{(n)}; (\tilde{h}, 0)\right)$ and $p_{\text{MD}}^{(n)}\left(R_{\boldsymbol{x}}^{(n)}; (\tilde{h}, h)\right)$, respectively. By [2, Theorem 13.4.1] there exists a further sub-sequence $\{n_{j_m}\}$ and a discriminant rule, $R_{\boldsymbol{x}}$, in the model which consists of a single observation from

$$\left\{\mathcal{N}_{2k}\left(\Lambda(\tilde{h}, h), \Lambda\right) \middle| (\tilde{h}, h) \in \mathbb{R}^{2k}\right\}, \tag{65}$$

such that for every $(\tilde{h}, h)$, where $h \neq 0$

$$p_{\text{MD}}^{(n_{j_m})}\left(R_{\boldsymbol{x}}^{(n_{j_m})}; (\tilde{h}, h)\right) \longrightarrow p_{\text{MD}}\left(R_{\boldsymbol{x}}; (\tilde{h}, h)\right) \tag{66}$$

and for $(\tilde{h}, 0)$

$$p_{\text{FA}}^{(n_{j_m})}\left(R_{\boldsymbol{x}}^{(n_{j_m})}; (\tilde{h}, 0)\right) \longrightarrow p_{\text{FA}}\left(R_{\boldsymbol{x}}; (\tilde{h}, 0)\right) \tag{67}$$

where the false-alarm and missed detection probabilities of $R_{\boldsymbol{x}}$ in (66) and (67) are defined with respect to the HT problem: $H_0 : \mathcal{N}_{2k}\left(\Lambda(\tilde{h}, 0), \Lambda\right)$ vs. $H_1 : \mathcal{N}_{2k}\left(\Lambda(\tilde{h}, h), \Lambda\right)$. The model (65) can be written as

$$\left\{\mathcal{N}_k\left(\rho\tilde{\mu}, \rho I\right) \times \mathcal{N}_k\left(\mu, I\right) \middle| (\tilde{\mu}, \mu) \in \mathbb{R}^{2k}\right\}, \tag{68}$$

where

$$\tilde{\mu} = J_{\theta_0}^{1/2} \tilde{h}.$$

Therefore, there exists a discriminant rule in the model (68) such that the false-alarm and missed-detection probabilities of $R_{\boldsymbol{x}}^{(n_{jm})}$ converges to those of that discriminant rule. For simplicity of notation we denote that discriminant rule by $R_{\boldsymbol{x}}$ as well. That is, for every $(\tilde{\mu}, \mu)$ where $\mu \neq 0$

$$p_{\text{MD}}^{(n_{jm})}\left(R_{\boldsymbol{x}}^{(n_{jm})}; (\tilde{\mu}, \mu)\right) \longrightarrow p_{\text{MD}}\left(R_{\boldsymbol{x}}; (\tilde{\mu}, \mu)\right) \tag{69}$$

and for $(\tilde{\mu}, 0)$

$$p_{\text{FA}}^{(n_{jm})}\left(R_{\boldsymbol{x}}^{(n_{jm})}; (\tilde{\mu}, 0)\right) \longrightarrow p_{\text{FA}}\left(R_{\boldsymbol{x}}; (\tilde{\mu}, 0)\right). \tag{70}$$

Note that for $\tilde{\mu}$ such that $\|\tilde{\mu}\| \geq d$, (70) ensures that $p_{\text{FA}}(R_{\boldsymbol{x}}; (\tilde{\mu}, 0)) \leq p_{\text{FA}}$, since $(\tilde{\mu}, 0)$ in (64) corresponds to the problem (20) and by assumption $R_{\boldsymbol{x}}^{(n)}$ has a guaranteed significance level $p_{\text{FA}}$. To complete the proof we will show that the missed-detection probability of the discriminant rule satisfying (69) and (70) in the equivalent model (which we will also denote by $R_{\boldsymbol{x}}$) cannot be less than $\underline{p}_{\text{MD}}^{\text{umm}}(p_{\text{FA}}, d, \rho, k)$. Since the model (68) is equivalent to the model

$$\left\{\mathcal{N}_k\left(\tilde{\mu}, \frac{1}{\rho}I\right) \times \mathcal{N}_k(\mu, I) \middle| (\tilde{\mu}, \mu) \in \mathbb{R}^{2k}\right\}, \tag{71}$$

there exists a discriminant rule $R_{\boldsymbol{x}}$ in the model (71) such that (69) and (70) hold. When $(\tilde{\mu}, \mu) = (a, a)$ or $(\tilde{\mu}, \mu) = (a, 0)$ for some $0 \neq a \in \mathbb{R}^k$, the limit model (71) corresponds to the problem in (13). For $a$ satisfying $d \leq \Delta$, where $\Delta = \|a\|$, notice that $R_{\boldsymbol{x}}$ is a discriminant rule with guaranteed significance level $p_{\text{FA}}$ in the normal location problem (13) with $\mu_1 = a$. Suppose to the contrary that

$$p < \underline{p}_{\text{MD}}^{\text{umm}}(p_{\text{FA}}, \Delta, \rho, k).$$

Then the missed-detection probability of $R_{\boldsymbol{x}}$ is better than $\underline{p}_{\text{MD}}^{\text{umm}}(p_{\text{FA}}, \Delta, \rho, k)$ for all $a$ such that $d \leq \Delta$. But this contradicts Theorem 1, therefore we must have

$$p \geq \underline{p}_{\text{MD}}^{\text{umm}}(p_{\text{FA}}, d, \rho, k),$$

for all $a$ such that $d \leq \Delta$.

∎

Finally, to complete the proof of Lemma 4 we have the following lemma which shows that two independent observations from a LAN model also form a LAN model.

*Lemma 7:* Let $\left(\Pi_\theta^{(n)} \big| \theta \in \Theta\right)_{n \in \mathbb{N}}$ be LAN at $\theta_0$ with norming matrices $r_n$ and information matrix $J_{\theta_0}$ and let

$$(\mathbf{W}, \mathbf{Y}) \sim \Pi_{\tilde{\theta}}^{(m)} \times \Pi_\theta^{(n)}.$$

If as $m(n), n \to \infty$ the matrices $[r_m r_n^{-1}]$ have a finite limit $\Gamma < \infty$, then the model corresponding to $(\mathbf{W}, \mathbf{Y})$ is LAN at $(\theta_0, \theta_0) \in \mathbf{R}^{2k}$ with norming block-matrices $\Psi_n$ and with information matrix $\Lambda$ which are given by the following $[2k \times 2k]$ block-matrices

$$\Psi_n = \begin{bmatrix} r_{n_m} & 0 \\ 0 & r_n \end{bmatrix}, \Lambda = \begin{bmatrix} \Gamma^t J_{\theta_0} \Gamma & 0 \\ 0 & J_{\theta_0} \end{bmatrix}.$$

*Proof:* We parametrize locally around $\theta_0$ by: $\theta_n = \theta_0 + r_n^{-1} h_n$ and $\tilde{\theta}_n = \theta_0 + r_n^{-1} \tilde{h}_n$. If $h_n \to h$ and $\tilde{h}_n \to \tilde{h}$ then using the Definition 13, the log likelihood-ratio is given by

$$\log \frac{L_{(\tilde{\theta}_n, \theta_n)}(\mathbf{W}, \mathbf{Y})}{L_{(\theta_0, \theta_0)}(\mathbf{W}, \mathbf{Y})} = h^t u_{n, \theta_0} - \frac{1}{2} h^t J_{\theta_0} h + \delta_n + \log \frac{L_{\tilde{\theta}_n}(\mathbf{W})}{L_{\theta_0}(\mathbf{W})}.$$

Notice that $r_n^{-1} \tilde{h}_n = r_m^{-1} q$, where $q = \left[r_m r_n^{-1} \tilde{h}_n\right]$ which under our assumption satisfies $q \to \Gamma \tilde{h}$. Therefore

$$\log \frac{L_{\tilde{\theta}_n}(\mathbf{W})}{L_{\theta_0}(\mathbf{W})} = (\Gamma \tilde{h})^t u_{m, \theta_0} - \frac{1}{2} (\Gamma \tilde{h})^t J_{\theta_0} (\Gamma \tilde{h}) + \tilde{\delta}_m$$

$$= \tilde{h}^t \left(\Gamma^t u_{m, \theta_0}\right) - \frac{1}{2} \tilde{h}^t \left(\Gamma^t J_{\theta_0} \Gamma\right) \tilde{h} + \tilde{\delta}_m.$$

Therefore the model corresponding to $(\mathbf{W}, \mathbf{Y})$ is LAN at $(\theta_0, \theta_0)$ with norming block-matrices $\Psi_n$ and with information matrix $\Lambda$.

∎

## APPENDIX C

### PROOF OF THEOREM 3

We prove Theorem 3 for the case where $\rho_m > 0$ for all $m$. The proof for the case where $\{\rho_m\}$ may have zero elements can be proved by continuity. For readability, the proof of Theorem 3 is split to three lemmas. The first lemma states $\mathcal{E}_m$ in terms of the non-centrality parameters $\vartheta_0, \vartheta_1$ in Theorem 1. The following two lemmas evaluate the missed-detection probability using normal approximation to the non-central $\mathcal{X}^2$-distribution.

*Lemma 8:* $\mathcal{E}_m$ in (37) satisfies

$$\mathcal{E}_m = E\left[W(\vartheta_0, \vartheta_1)\right] + O\left(\frac{1}{\sqrt{E(b_m)}}\right),$$

where

$$W(\vartheta_0, \vartheta_1) = \frac{\vartheta_1 - \vartheta_0}{\sqrt{b_m}},$$

$$b_m = 2k_m + 4\vartheta_1.$$

*Lemma 9:* The missed-detection probability satisfies

$$1 - p_{\mathrm{MD}}^{(m)} = E\left[Q\left(Q^{-1}\left(p_{\mathrm{FA}}^{(m)}\right) - \mathcal{E}_m\right)\right] + O\left(\frac{1}{\sqrt{\max\{k_m, \rho_m\}}}\right), \tag{72}$$

where

$$p_{\mathrm{FA}}^{(m)} = Q\left(\sqrt{1-U}Q^{-1}(\tilde{p}_{\mathrm{FA}})\right),$$

$$\tilde{p}_{\mathrm{FA}} = p_{\mathrm{FA}} + O\left(\frac{1}{\sqrt{\max\{k_m, \vartheta_0\}}}\right), \tag{73}$$

$$U = \frac{4W(\vartheta_0, \vartheta_1)}{\sqrt{b_m}}. \tag{74}$$

*Lemma 10:* The mean on the RHS of (72) satisfies

$$E\left[Q\left(Q^{-1}\left(p_{\mathrm{FA}}^{(m)}\right) - \mathcal{E}_m\right)\right] = Q\left(Q^{-1}\left(p_{\mathrm{FA}} + \tilde{\psi}\right) - \mathcal{E}_m\right) + O\left(\frac{1}{\sqrt{\max\{k_m, \rho_m\}}}\right),$$

where

$$\tilde{\psi} = O\left(\frac{1}{\sqrt{\max\{k_m, \rho_m\}}}\right).$$

## A. Proof of Lemma 8

From (48), (49) and by the properties of the $\mathcal{X}^2$-distribution we can obtain the means and variances of $\vartheta_0, \vartheta_1$. Specifically,

$$\begin{aligned} E(\vartheta_1) &= \rho_m(k_m + \lambda_1) \\ &= \rho_m k_m + (1 + \rho_m)^2 \Delta_m^2, \\ \mathrm{Var}(\vartheta_1) &= \rho_m^2(2k_m + 4\lambda_1) \\ &= 2\rho_m^2 k_m + 4\rho_m(1+\rho_m)^2 \Delta_m^2 \end{aligned} \tag{75}$$

and

$$E(\vartheta_0) = \rho_m(k_m + \lambda_0)$$

$$=\rho_m k_m + \rho_m^2 \Delta_m^2,$$

$$\text{Var}(\vartheta_0) = \rho_m^2 (2k_m + 4\lambda_0)$$

$$= 2\rho_m^2 k_m + 4\rho_m^3 \Delta_m^2$$

$$\leq \text{Var}(\vartheta_1). \tag{76}$$

The difference between the non-centrality parameters satisfies

$$\vartheta_1 - \vartheta_0 = \Delta_m^2 + 2\rho_m X^t \mu_1.$$

Notice that $X^t \mu_1 \sim \mathcal{N}\left(\Delta_m^2, \frac{\Delta_m^2}{\rho_m}\right)$. Therefore

$$\vartheta_1 - \vartheta_0 \sim \mathcal{N}\left(\Delta_m^2(1 + 2\rho_m), 4\rho_m \Delta_m^2\right). \tag{77}$$

Using (75) and (77) it follows that $\mathcal{E}_m$ satisfies

$$\mathcal{E}_m = \frac{E(\vartheta_1 - \vartheta_0)}{\sqrt{E(b_m)}}.$$

We shall use the second order Taylor approximation of $W(\vartheta_0, \vartheta_1)$. The partial derivatives of $W$ will be denoted by the corresponding subscripts and are given by:

$$W_{\vartheta_0} = -\frac{1}{\sqrt{b_m}}$$

$$W_{\vartheta_1} = \frac{1}{\sqrt{b_m}} - \frac{2(\vartheta_1 - \vartheta_0)}{b_m^{3/2}}$$

$$= \frac{1}{\sqrt{b_m}} - \frac{2W}{b_m}.$$

The second order derivatives are given by

$$W_{\vartheta_0,\vartheta_0} = 0$$

$$W_{\vartheta_0,\vartheta_1} = \frac{2}{b_m^{3/2}}$$

$$W_{\vartheta_1,\vartheta_1} = -\frac{2}{b_m^{3/2}} - 2\left[\frac{1}{b_m^{3/2}} - \frac{6W}{b_m^2}\right]$$

$$= -\frac{4}{b_m^{3/2}} + \frac{12W}{b_m^2}.$$

The second order Taylor approximation about $(E(\vartheta_0), E(\vartheta_1))$ yields

$$W(\vartheta_0, \vartheta_1) = \mathcal{E}_m - \frac{\vartheta_0 - E(\vartheta_0)}{\sqrt{E(b_m)}} + \left[\frac{1}{\sqrt{E(b_m)}} - \frac{2\mathcal{E}_m}{E(b_m)}\right](\vartheta_1 - E(\vartheta_1)) + \delta, \tag{78}$$

where

$$\delta = O\left(\frac{(\vartheta_1 - E(\vartheta_1))^2}{(E(b_m))^{3/2}} + \frac{(\vartheta_1 - E(\vartheta_1))(\vartheta_0 - E(\vartheta_0))}{(E(b_m))^{3/2}}\right).$$

Taking the expectation of (78) yields

$$E[W(\vartheta_0, \vartheta_1)] = \mathcal{E}_m + E(\delta).$$

To complete the proof we show that

$$E(\delta) = O\left(\frac{1}{\sqrt{E(b_m)}}\right). \tag{79}$$

Using the Cauchy-Schwarz inequality and the fact that $Var(\vartheta_1) = E(b_m) - \rho_m k_m$, we get

$$E\left[\frac{(\vartheta_1 - E(\vartheta_1))^2}{(E(b_m))^{3/2}}\right] = \frac{E(b_m) - \rho_m k_m}{(E(b_m))^{3/2}}$$

$$= O\left(\frac{1}{\sqrt{E(b_m)}}\right), \tag{80}$$

$$\frac{\text{Cov}(\vartheta_0, \vartheta_1)}{(E(b_m))^{3/2}} \leq \frac{\sqrt{\text{Var}(\vartheta_0)\text{Var}(\vartheta_1)}}{(E(b_m))^{3/2}}$$

$$\leq \frac{\text{Var}(\vartheta_1)}{(E(b_m))^{3/2}} \tag{81}$$

$$= O\left(\frac{1}{\sqrt{E(b_m)}}\right). \tag{82}$$

where the inequality (81) is due to (76). Clearly (80) and (82) yield (79).

*B. Proof of Lemma 9*

For any given $x$ let

$$\tilde{p}_{\text{FA}} = Q\left(\frac{Q_{(k_m),\vartheta_0}^{-1}(p_{\text{FA}}) - k_m - \vartheta_0}{\sqrt{2k_m + 4\vartheta_0}}\right).$$

Then for

$$T = k_m + \vartheta_0 + \sqrt{2k_m + 4\vartheta_0}Q^{-1}(\tilde{p}_{\text{FA}}),$$

we have

$$Q_{(k_m),\vartheta_0}(T) = p_{\text{FA}}.$$

We use the normal approximation for the non-central $\mathcal{X}^2$-distribution. Specifically,

$$Q_{(k_m),\vartheta_0}\left(\sqrt{2k_m + 4\vartheta_0}Q^{-1}(p_{\text{FA}}) + k_m + \vartheta_0\right) = p_{\text{FA}} + O\left(\frac{1}{\sqrt{\max\{k_m, \vartheta_0\}}}\right). \quad (83)$$

Indeed the correction term may be smaller than the $O(1/\sqrt{k_m})$ guaranteed by the Berry-Esseen Theorem since the correction term is also $O\left(1/\sqrt{\vartheta_0}\right)$, see [30, Page 466]. It follows that

$$Q^{-1}_{(k_m),\vartheta_0}(p_{\text{FA}}) = T, \quad (84)$$

$$\tilde{p}_{\text{FA}} = p_{\text{FA}} + O\left(\frac{1}{\sqrt{\max\{k_m, \vartheta_0\}}}\right).$$

Using (84) and approximations similar to (83) for the conditional missed-detection probability yields

$$\begin{aligned}
1 - p^{(m)}_{\text{MD}}(\underline{R}_{\boldsymbol{x}} \mid x; \mu_1) &= Q_{(k_m),\vartheta_1}\left(Q^{-1}_{(k_m),\vartheta_0}(p_{\text{FA}})\right) \\
&= Q\left(\frac{Q^{-1}_{(k_m),\vartheta_0}(p_{\text{FA}}) - k_m - \vartheta_1}{\sqrt{b_m}}\right) + O\left(\frac{1}{\sqrt{\max\{k_m, \vartheta_1\}}}\right) \\
&= Q\left(\frac{\sqrt{2k_m + 4\vartheta_0}Q^{-1}(\tilde{p}_{\text{FA}})}{\sqrt{b_m}} - W(\vartheta_0, \vartheta_1)\right) + O\left(\frac{1}{\sqrt{\max\{k_m, \vartheta_1\}}}\right) \\
&= Q\left(\sqrt{1-U}Q^{-1}(\tilde{p}_{\text{FA}}) - W(\vartheta_0, \vartheta_1)\right) + O\left(\frac{1}{\sqrt{\max\{k_m, \vartheta_1\}}}\right), \quad (85)
\end{aligned}$$

We rewrite (85) as

$$1 - p^{(m)}_{\text{MD}}(\underline{R}_{\boldsymbol{x}} \mid x; \mu_1) + \xi_1 = Q\left(Q^{-1}\left(p^{(m)}_{\text{FA}}\right) - W(\vartheta_0, \vartheta_1)\right), \quad (86)$$

where

$$\xi_1 = O\left(\frac{1}{\sqrt{\max\{k_m, \vartheta_1\}}}\right).$$

Using the Taylor approximation of the RHS of (86) about $\left(Q^{-1}\left(p^{(m)}_{\text{FA}}\right) - \mathcal{E}_m\right)$ yields

$$Q\left(Q^{-1}\left(p^{(m)}_{\text{FA}}\right) - W(\vartheta_0, \vartheta_1)\right) = Q\left(Q^{-1}\left(p^{(m)}_{\text{FA}}\right) - \mathcal{E}_m\right) + \xi_2, \quad (87)$$

where

$$\xi_2 = -\phi(c)\left[\mathcal{E}_m - W(\vartheta_0, \vartheta_1)\right],$$

and $c$ lies between $\left(Q^{-1}\left(p_{\text{FA}}^{(m)}\right) - \mathcal{E}_m\right)$ and $\left(Q^{-1}\left(p_{\text{FA}}^{(m)}\right) - W(\vartheta_0, \vartheta_1)\right)$. From (86) and (87) we get

$$1 - p_{\text{MD}}^{(m)} + E(\xi_1) = E\left[Q\left(Q^{-1}\left(p_{\text{FA}}^{(m)}\right) - \mathcal{E}_m\right)\right] + E(\xi_2). \tag{88}$$

According to Lemma 8

$$E(\xi_2) = O\left(\frac{1}{\sqrt{E(b_m)}}\right). \tag{89}$$

Now we quantify the other terms in (88) separately starting with $E(\xi_1)$. To this end we construct a confidence interval for $\vartheta_1/\rho_m$ using the normal approximation to the non-central $\mathcal{X}^2$-distribution. Specifically, let

$$G_1 = \left[k_m + \lambda_1 \pm s_1\sqrt{2k_m + 4\lambda_1}\right],$$

$$s_1 = Q^{-1}\left(\frac{1}{\sqrt{k_m + \lambda_1}}\right).$$

By (75) the confidence level of $G_1$ satisfies

$$P\left(\frac{\vartheta_1}{\rho_m} \notin G_1\right) = 2Q(s_1) + O\left(\frac{1}{\sqrt{\max\{k_m, \lambda_1\}}}\right) \tag{90}$$

$$= O\left(\frac{1}{\sqrt{\max\{k_m, \lambda_1\}}}\right).$$

Furthermore, the event $A_1 = \{\vartheta_1/\rho_m \in G_1\}$ implies

$$\vartheta_1 = \rho_m(k_m + \lambda_1) + O\left(\rho_m s_1 \sqrt{k_m + \lambda_1}\right). \tag{91}$$

By assumption $\{\mathcal{E}_m\}$ is bounded away from $0$ and $\infty$. This implies that

$$O\left(\frac{1}{\sqrt{k_m + (1+\rho_m)\Delta_m^2}}\right) = O\left(\frac{1}{\sqrt{k_m + \rho_m \Delta_m^2}}\right) \tag{92}$$

$$= O\left(\frac{1}{\sqrt{\max\{k_m, \rho_m\}}}\right).$$

To see this write

$$\mathcal{E}_m = \frac{\Delta_m^2 \sqrt{1 + 2\rho_m}}{\sqrt{k_m + \frac{4(1+\rho_m)^2}{(1+2\rho_m)}\Delta_m^2}}.$$

Note that if

$$k_m \ll (1+\rho_m)\Delta_m^2, \tag{93}$$

then $\mathcal{E}_m = \Theta(\Delta_m)$ which is bounded away from $\infty$, thus clearly (92) holds. On the other hand, if (93) does not hold then it must be that

$$\Delta_m^2 = \Theta\left(\sqrt{\frac{k_m}{1+\rho_m}}\right),$$

which yields (92). Using (90), (91) and (92) yields

$$\begin{aligned}
E(\xi_1) &= O\left(\frac{1}{\sqrt{b_m}}\right) P(A_1) + O\left(\frac{1}{\sqrt{k_m}}\right) P(A_1^c) \\
&= O\left(\frac{1}{\sqrt{b_m}}\right) + O\left(\frac{1}{\sqrt{\max\{k_m^2, k_m\lambda_1\}}}\right) \\
&= O\left(\frac{1}{\sqrt{k_m + (1+\rho_m)\Delta_m^2}}\right) \\
&= O\left(\frac{1}{\sqrt{\max\{k_m, \rho_m\}}}\right).
\end{aligned} \quad (94)$$

Clearly (88), (89) and (94) yield (72).

*C. Proof of Lemma 10*

We construct a confidence interval for $(\vartheta_1 - \vartheta_0)$ using (77) as follows

$$\tilde{G} = \left[\Delta_m^2(1+2\rho_m) \pm \tilde{s}\sqrt{4\rho_m\Delta_m^2}\right],$$

$$\tilde{s} = Q^{-1}\left(\frac{1}{\sqrt{\Delta_m^2(1+2\rho_m)}}\right).$$

The event $\tilde{A} = \{(\vartheta_1 - \vartheta_0) \in \tilde{G}\}$ implies

$$\vartheta_1 - \vartheta_0 = \Delta_m^2(1+2\rho_m) + O\left(\tilde{s}\sqrt{\rho_m\Delta_m^2}\right). \quad (95)$$

The confidence level of $\tilde{G}$ satisfies

$$\begin{aligned}
P\left((\vartheta_1 - \vartheta_0) \notin \tilde{G}\right) &= 2Q(\tilde{s}) \\
&= O\left(\frac{1}{\sqrt{\Delta_m^2(1+2\rho_m)}}\right).
\end{aligned}$$

According to (91) and (95) we can quantify $W$ as follows

$$1_{\{\tilde{A} \cap A_1\}} \cdot W = \frac{\Delta_m^2(1+2\rho_m) + O\left(\tilde{s}\sqrt{\rho_m\Delta_m^2}\right)}{\sqrt{E(b_m) + O\left(\rho_m s_1 \sqrt{k_m + \lambda_1}\right)}} \quad (96)$$

$$=O\left(\mathcal{E}_m\right)$$
$$=O\left(1\right).$$

From (74) and (96), it follows that $U$ satisfies

$$1_{\{\tilde{A} \cap A_1\}} \cdot U = O\left(\frac{1}{\sqrt{E(b_m)}}\right). \tag{97}$$

Finally, we construct a confidence interval for $\vartheta_0/\rho_m$ as follows

$$G_0 = \left[k_m + \lambda_0 \pm s_0\sqrt{2k_m + 4\lambda_0}\right],$$
$$s_0 = Q^{-1}\left(\frac{1}{\sqrt{k_m + \lambda_0}}\right).$$

The confidence level of $G_0$ satisfies

$$P\left(\frac{\vartheta_0}{\rho_m} \notin G_0\right) = 2Q(s_0) + O\left(\frac{1}{\sqrt{\max\{k_m, \rho_m\Delta_m^2\}}}\right)$$
$$= O\left(\frac{1}{\sqrt{k_m + \rho_m\Delta_m^2}}\right).$$

The event $A_0 = \{\vartheta_0/\rho_m \in G_0\}$ implies

$$\vartheta_0 = \rho_m(k_m + \lambda_0) + O\left(\rho_m s_0 \sqrt{k_m + \lambda_0}\right). \tag{98}$$

Using the confidence intervals, we will show that $p_{\text{FA}}^{(m)}$ satisfies

$$1_{\{A_0 \cap A_1 \cap \tilde{A}\}} \cdot p_{\text{FA}}^{(m)} = 1_{\{A_0 \cap A_1 \cap \tilde{A}\}} \cdot p_{\text{FA}} + O\left(\frac{1}{\sqrt{\max\{k_m, \rho_m\Delta_m^2\}}}\right). \tag{99}$$

To see this denote

$$g(u) = Q\left(\sqrt{1-u}Q^{-1}(\tilde{p}_{\text{FA}})\right).$$

Using the Taylor approximation of $g(u)$ about 0 yields

$$g(u) = \tilde{p}_{\text{FA}} + \gamma,$$

where

$$\gamma = \phi\left(\sqrt{1-c'}Q^{-1}(\tilde{p}_{\text{FA}})\right) \cdot \frac{uQ^{-1}(\tilde{p}_{\text{FA}})}{2\sqrt{1-c'}},$$

and where $c'$ lies between 0 and $u$. According to (97)

$$\gamma = O\left(\frac{\phi\left(Q^{-1}(\tilde{p}_{\text{FA}})\right)Q^{-1}(\tilde{p}_{\text{FA}})}{\sqrt{E(b_m)}}\right).$$

Now, since for any $\tilde{p}_{\text{FA}}$

$$\left|\phi\left(Q^{-1}(\tilde{p}_{\text{FA}})\right) Q^{-1}(\tilde{p}_{\text{FA}})\right| \leq 1,$$

we conclude that

$$\gamma = O\left(\frac{1}{\sqrt{E(b_m)}}\right). \tag{100}$$

But notice that (73), (98) and (100) imply that

$$\begin{aligned}
1_{\{A_0 \cap A_1 \cap \tilde{A}\}} \cdot p_{\text{FA}}^{(m)} &= 1_{\{A_0 \cap A_1 \cap \tilde{A}\}} \cdot \tilde{p}_{\text{FA}} + \gamma \\
&= 1_{\{A_0 \cap A_1 \cap \tilde{A}\}} \cdot p_{\text{FA}} + O\left(\frac{1}{\sqrt{\max\{k_m, \rho_m \Delta_m^2\}}}\right) + \gamma \\
&= 1_{\{A_0 \cap A_1 \cap \tilde{A}\}} \cdot p_{\text{FA}} + O\left(\frac{1}{\sqrt{\max\{k_m, \rho_m \Delta_m^2\}}}\right),
\end{aligned} \tag{101}$$

thus (99) holds. Now, the probability of the intersection of events $A_0 \cap A_1 \cap \tilde{A}$ satisfies

$$P\left(A_0 \cap A_1 \cap \tilde{A}\right) = 1 + O\left(\frac{1}{\sqrt{\max\{k_m, \lambda_0\}}}\right) + O\left(\frac{1}{\sqrt{\max\{k_m, \lambda_1\}}}\right) + O\left(\frac{1}{\sqrt{\Delta_m^2(1 + 2\rho_m)}}\right)$$

$$= 1 + O\left(\frac{1}{\sqrt{\max\{k_m, \rho_m \Delta_m^2\}}}\right). \tag{102}$$

Using (92), (101) and (102) it follows that

$$E\left[Q\left(Q^{-1}\left(p_{\text{FA}}^{(m)}\right) - \mathcal{E}_m\right)\right] = 1_{\{A_0 \cap A_1 \cap \tilde{A}\}} \cdot E\left[Q\left(Q^{-1}\left(p_{\text{FA}}^{(m)}\right) - \mathcal{E}_m\right)\right] + O\left(\frac{1}{\sqrt{\max\{k_m, \rho_m \Delta_m^2\}}}\right)$$

$$= Q\left(Q^{-1}\left(p_{\text{FA}} + \tilde{\psi}\right) - \mathcal{E}_m\right) + O\left(\frac{1}{\sqrt{\max\{k_m, \rho_m\}}}\right),$$

where

$$\tilde{\psi} = O\left(\frac{1}{\sqrt{\max\{k_m, \rho_m \Delta_m^2\}}}\right)$$

$$= O\left(\frac{1}{\sqrt{\max\{k_m, \rho_m\}}}\right)$$

as required.

# REFERENCES

[1] T. M. Cover and J. A. Thomas. *Elements of information theory*. Wiley, New York, 1991.

[2] E. L Lehmann and J. P Romano. *Testing statistical hypotheses*. Springer, third edition, 2005.

[3] M. Feder and N. Merhav. Universal composite hypothesis testing: a competitive minimax approach. *IEEE Trans. Inform. Theory*, 48(6):1504–1517, Jun 2002.

[4] Y. Polyanskiy, H. V. Poor, and S. Verdú. Channel coding rate in the finite blocklength regime. *IEEE Trans. Inform. Theory*, 56(5):2307–2359, May 2010.

[5] O. Zeitouni, J. Ziv, and N. Merhav. When is the generalized likelihood ratio test optimal? *IEEE Trans. Inform. Theory*, 38(5):1597–1602, Sep 1992.

[6] E. Levitan and N. Merhav. A competitive Neyman-Pearson approach to universal hypothesis testing with applications. *IEEE Trans. Inform. Theory*, 48(8):2215–2229, Aug 2002.

[7] M. Gutman. Asymptotically optimal classification for multiple tests with empirically observed statistics. *IEEE Trans. Inform. Theory*, 35(2):401–408, March 1989.

[8] V. Strassen. Asymptotische abschatzungen in Shannon's informationstheorie. *Trans. Third Prauge Conf. Information theory, Czechoslovak Academy of Sciences*, pages 689–723, 1962.

[9] I. Sason. Moderate deviations analysis of binary hypothesis testing. In *2012 IEEE International Symposium on Information Theory*, pages 821–825, July 2012.

[10] L. Zhou, V. Y. F. Tan, and M. Motani. Second-order asymptotically optimal statistical classification. *Information and Inference: A Journal of the IMA*, 01 2019.

[11] Y. Huang and P. Moulin. Strong large deviations for composite hypothesis testing. In *2014 IEEE International Symposium on Information Theory*, pages 556–560, June 2014.

[12] Y. Li, S. Nitinawarat, and V. V. Veeravalli. Universal outlier hypothesis testing. *IEEE Trans. Inform. Theory*, 60(7):4066–4082, July 2014.

[13] B. Mauro and T Benedetta. Binary hypothesis testing game with training data. *IEEE Trans. Inform. Theory*, 60(8):4848–4866, Aug 2014.

[14] F. Farnia and D. Tse. A minimax approach to supervised learning. In D. D. Lee, M. Sugiyama, U. V. Luxburg, I. Guyon, and R. Garnett, editors, *Advances in Neural Information Processing Systems 29*, pages 4240–4248. Curran Associates, Inc., 2016.

[15] S.-L. Huang, A. Makur, G.W. Wornell, and L. Zheng. On Universal Features for High-Dimensional Learning and Inference. *ArXiv:1911.09105*, Nov. 2019.

[16] Y. I. Ingster and I. A. Suslina. *Nonparametrie goodness-of-fit testing under Gaussian models*. Springer, 2003.

[17] A. W. van der Vaart. *Asymptotic statistics*. Cambridge University Press, Cambridge, UK, 1998.

[18] J. C. Berry. Minimax estimation of a bounded mean vector. *Journal of Multivariate Analysis*, 35(1):130–139, October 1990.

[19] K.V. Mardia, J.T. Kent, and J.M. Bibby. *Multivariate analysis*. Academic Press, 1979.

[20] G.G. Roussas. *Contiguity of Probability Measures.*. Cambridge University Press, 1972.

[21] G.G. Roussas. Asymptotic distribution of the log-likelihood function for stochastic processes. *Zeitschrift für Wahrscheinlichkeitstheorie und Verwandte Gebiete*, 47:31–46, 1979.

[22] M. Hallin, D. M. Mason, D. Pfeifer, and J. G. Steinebach. *Mathematical Statistics and Limit Theorems*. Springer, 2015.

[23] A. Brouste and M. Fukasawa. Local asymptotic normality property for fractional Gaussian noise under high-frequency observations. *Annals of Statistics*, 46(5):2045–2061, 2018.

[24] I.A. Ibragimov and R.Z. Has'minskii. *Statistical estimation asymptotic theory*. Springer, 1981.

[25] A.W. van der Vaart. The statistical work of Lucien Le Cam. *The Annals of Statistics*, 30(3):631 – 682, 2002.

[26] I. Csiszár and P. C. Shields. *Information theory and statistics: A tutorial*. Now Publishers Inc, 2004.

[27] T.S. Ferguson. An inconsistent maximum likelihood estimate. *Journal of the American Statistical Association*, 77:831–834, 1982.

[28] H. Cramér. *Mathematical Methods of Statistics*. Princeton University Press, 1946.

[29] G. A. Rempała and J. Wesołowski. Double asymptotics for the chi-square statistic. *ArXiv e-prints*, September 2016.

[30] N.L. Johnson, S.N. Kotz, and N. Balakrishnan. *Continuous Univariate Distributions*, volume 2. Wiley, second edition, 1995.